\documentclass[11pt]{article}
\usepackage[table]{xcolor}
\usepackage{graphicx}
\usepackage{amsmath,amsfonts,amssymb}
\usepackage{graphicx}
\usepackage{multirow}
\usepackage{tikz,pgf}
\usepackage{url}
\usepackage{amsmath}
\usepackage{graphicx}
\usepackage{epsfig}
\usepackage{epstopdf}
\usepackage{color}
\usepackage{enumitem}
\def\tto{\;{\lower 1pt \hbox{$\rightarrow$}}\kern -10pt
\hbox{\raise 2pt \hbox{$\rightarrow$}}\;}

\def\Bar{\overline}
\def\ra{\rangle}
\def\la{\langle}

\def\epsilon{\varepsilon}

\def\h{\hfill\Box}
\def\R{\Bbb R}

\def\ox{\bar{x}}

\def\oy{\bar{y}}
\def\oz{\bar{z}}

\def\ou{\bar{u}}

\def\ri{\mbox{\rm ri}\,}

\def\gph{\mbox{\rm gph}\,}

\def\epi{\mbox{\rm epi}\,}

\def\dom{\mbox{\rm dom}\,}

\def\cone{\mbox{\rm cone}}
\def\rge{\mbox{\rm rge}\,}

\def\sqri{\mbox{\rm sqri}\,}

\def\h{\hfill\square}

\def\oR{\Bar{\R}}
\def\s{\square}

\def\oR{\Bar{\R}}

\DeclareMathOperator{\qri}{qri}
\setlist[enumerate,1]{itemsep=0.0ex,parsep=0.5ex,label={\rm(\alph*)},leftmargin=*, align=left}
\newcounter{lk}

\definecolor{gls}{rgb}{0.0, 0.5, 1.0}
\begin{document}
\begin{center}
{\sc\bf  Fenchel Conjugate of Set-Valued Mappings}\\[1ex]
{\sc Nguyen Mau   Nam}\footnote{Fariborz Maseeh Department of Mathematics and Statistics, Portland State University, Portland, OR
97207, USA (mnn3@pdx.edu). Research of this author was partly supported by the USA National Science Foundation under grant DMS-2136228.},
{\sc Gary Sandine}\footnote{Fariborz Maseeh Department of
 Mathematics and Statistics, Portland State University, Portland, OR 97207, USA (gsandine@pdx.edu).},
{\sc Nguyen Nang Thieu}\footnote{Institute of Mathematics, Vietnam Academy of Science and Technology,
Hanoi, Vietnam \& The State University of New York - SUNY, Korea (nnthieu@math.ac.vn).},
{\sc Nguyen Dong Yen}\footnote{Institute of Mathematics, Vietnam Academy of Science and Technology, 18 Hoang Quoc Viet,
Hanoi 10307 (ndyen@math.ac.vn).  Research of this author was supported by Vietnam Academy of Science and Technology under project number NVCC01.07/22-23.}
\\[2ex]
 {\bf Dedicated to Prof. Boris Mordukhovich on the occasion of his 75th birthday}

\end{center}
\small{\bf Abstract.} In this paper, we present a novel concept of the Fenchel conjugate for set-valued mappings and investigate its properties in finite and infinite dimensions. After establishing the fundamental properties of the Fenchel conjugate for set-valued mappings, we derive its main calculus rules in various settings. Our approach is geometric and draws inspiration from the successful application of this method by  B.~S.~Mordukhovich and coauthors in variational and convex analysis. Subsequently, we demonstrate that our new findings for the Fenchel conjugate of set-valued mappings can be utilized to obtain many old and new calculus rules of convex generalized differentiation in both finite and infinite dimensions.  \\[1ex]

\newtheorem{theorem}{Theorem}[section]
\newtheorem{proposition}[theorem]{Proposition}
\newtheorem{remark}[theorem]{Remark}
\newtheorem{lemma}[theorem]{Lemma}
\newtheorem{corollary}[theorem]{Corollary}
\newtheorem{definition}[theorem]{Definition}
\newtheorem{example}[theorem]{Example}
\renewcommand{\theequation}{\thesection.\arabic{equation}}
\normalsize

\section{Introduction}
\setcounter{equation}{0}

The concept of the conjugate function was first introduced and extensively studied by W.~Fenchel in \cite{fenchel1949,fenchel1953} for convex functions on finite-dimensional spaces. After Fenchel's pioneering work in finite dimensions, A.~Br\o ndsted~\cite{brondsted1964}, J.~J.~Moreau~\cite{moreau1962,moreau1962a}, and R.~T.~Rockafellar~\cite{r1966} studied various properties of conjugate functions in infinite-dimensional spaces. This concept is now broadly known as the {\em Fenchel conjugate} for functions.   It plays a crucial role in the formulation of Fenchel dual problems in convex optimization, serving as a gateway to a broader topic called duality theory.

It turns out that the Fenchel conjugate for convex functions has close relationships with convex generalized differential theory. In particular, calculus rules for the Fenchel conjugate of convex functions can be used to  derive major calculus rules for subdifferentials of convex functions in both finite and infinite dimensions. Another important concept of modern convex analysis is the notion of {\em coderivative}, which was introduced by B.~Mordukhovich in the early 1980s for set-valued mappings without requiring convexity. Coderivative calculus rules have been developed for both convex and nonconvex cases in both finite and infinite dimensions, constituting a major part of convex and variational analysis with their applications to optimization; see~\cite{m-book,m-book1,bmn2022,rw} and the references therein.

The important role of the Fenchel conjugate of convex functions in convex optimization has motivated the search for a notion of the Fenchel conjugate of set-valued mappings with convex graphs. To the best of our knowledge, the first effort was made by W.~Song in his work~\cite{Song1998}, where he introduced the notion of the Fenchel conjugate mapping for set-valued mappings. Further studies  in this direction can be found in~\cite{bot1999}. This concept of the Fenchel conjugate of set-valued mappings was defined on partially ordered vector spaces, with the main objective of developing a duality theory in set-valued optimization. Furthermore, it appears that the Fenchel conjugate of set-valued mappings is not closely related to the coderivatives of the set-valued mappings, unlike the  relationships between the Fenchel conjugate and the subdifferential of convex functions.

In the present paper, we introduce a new notion of the Fenchel conjugate for set-valued mappings in the framework of topological vector spaces and investigate its properties. Based on a geometric approach which was successfully used by B.~S.~Mordukhovich and coauthors (see,  e.g.,~\cite{m-book1,m-book,bmn,bmn2022}), we build comprehensive calculus rules for this new notion for convex set-valued mappings in both finite and infinite dimensions. We show that  our notion of the Fenchel conjugate mapping has close relationships with the coderivative notion and also with the Fenchel conjugate of convex functions. Based on these relationships, we derive  refined calculus rules for coderivatives of convex set-valued mappings as well as  enhanced subdifferential and Fenchel conjugate calculus rules for convex  functions. Our investigations in this paper are based on some initial recent work in finite dimensions from~\cite{HNY2023}.

The paper is structured as follows. In Sect.~2, we provide basic definitions and some elements of convex analysis used throughout the paper. Section~3 focuses on fundamental properties of the Fenchel conjugate of set-valued mappings. Section~4 is devoted to developing major calculus rules for the Fenchel conjugate of convex set-valued mappings in finite dimensions. Comprehensive calculus rules for this notion are developed in  Sect.~5.  Then, in Sect.~6, we use our new results to derive old and new coderivative calculus rules for set-valued mappings as well as subdifferential and Fenchel conjugate calculus rules for convex functions. Some concluding remarks  are given in Sect.~\ref{s:Concluding}.

\section{Basic Definitions and Preliminaries}\label{Sect_2}
\setcounter{equation}{0}

Throughout the paper, we use notation and concepts of convex analysis from~\cite{bmn,bmn2022}. The topological dual of a real topological vector space $X$ is denoted by $X^*$. For $x^*\in X^*$ and $x\in X$, we use $\la x^*,x\ra$ for the function value $x^*(x)$.  In the sequel, unless otherwise stated, $X^*$ is equipped with the weak$^*$ topology (see~\cite[Definition~1.107 and Subsection~1.2.2]{bmn2022} for details). For $A\subset  X$, the cone generated by $A$ is the set $\cone A=\{ta\; |\; t\geq 0, a\in A\}$. The set $\oR=\R\cup \{\pm \infty\}$ is the extended real line. If $M$ is a matrix, then $M^T$ stands for the transpose of $M$. For simplicity of presentation, we use $X$, $Y$, and $Z$ to denote real locally convex Hausdorff topological vector spaces which are nontrivial.

We now recall some  basic definitions from convex analysis. As mentioned  above, for simplicity, we only consider sets, functions, and mappings defined on locally convex topological vector spaces. The reader can find more details  in~\cite{bmn,bmn2022,r,zalinescu2002}.

A  set $\Omega$ in $X$ is called a {\em convex set} if $\lambda x+(1-\lambda)u\in \Omega$ whenever $x, u \in \Omega$ and $0\leq \lambda\leq 1$.

Given a function $f\colon X\to \oR$, the {\em epigraph} and the {\em effective domain} of $f$ are defined respectively by
\begin{align*}
	& \epi f=\big\{(x, \lambda)\in X\times \R\; \big |\; f(x)\leq \lambda\big\}\ \; \mbox{\rm and }\;
	\dom f=\big\{x\in X\; \big |\; f(x)<\infty\big\}.
\end{align*}
We say that $f$ is a {\em convex function} if its epigraph is a convex set in $X\times \R$, and we say that $f$ is a {\em proper function} if $\dom f\neq\emptyset$ and $-\infty<f(x)$ for all $x\in X$. The function $f$ is said to be {\em lower semicontinuous/closed} if $\epi f$ is closed.

The {\em subdifferential} of a convex function $f\colon X\to \oR$ at $\ox$ with $f(\ox)\in \R$ is given by
\begin{equation*}
	\partial f(\ox)=\big\{x^*\in X^*\; \big |\; \la x^*, x-\ox\ra\leq f(x)-f(\ox)\ \; \mbox{\rm for all }\; x\in X\big\}.
\end{equation*}
Consider a convex function $h\colon X^*\to \oR$ with  $x^*\in X^*$ and $h(x^*)\in \R$. Then we define
\begin{equation*}
	\partial h(x^*)=\big\{x\in X\; \big |\; \la u^*-x^*, x\ra \leq h(u^*)-h(x^*)\ \; \mbox{\rm for all }\; u^*\in X^*\big\}.
\end{equation*}

Given a nonempty subset $\Omega\subset X$, the {\em support function} $\sigma_\Omega\colon X^*\to \oR$ of $\Omega$ is defined by
\begin{equation}\label{support function}		\sigma_\Omega(x^*)=\sup\big\{\langle x^*,x \rangle\mid x\in \Omega\big\},\ \; x^*\in X^*.
\end{equation}
It can be shown that the support function of a nonempty set $\Omega$ is a proper convex function
whether the set is convex or not.

Let $\Omega\subset X$ be a subset of $X$. The {\em indicator function} $\delta_\Omega\colon X\to \oR$ associated with $\Omega$ is defined by
\begin{equation*}
	\delta_\Omega(x)=
	\begin{cases}
		0\quad&\mbox{\rm if}\; x\in\Omega,\\
		\infty&\mbox{\rm otherwise.}
	\end{cases}
\end{equation*}
Given two proper functions $g, h\colon X\to \oR$, define their {\em infimal convolution} by
\begin{equation*}
	(g\s h)(x)=\inf\big\{g(x_1)+h(x_2)\; \big|\; x_1+x_2=x\big\}, \ \; x\in X.
\end{equation*}
For a multifunction/set-valued mapping $F\colon X\tto Y$, the {\em graph}, the {\em effective domain}, and the {\em range} of $F$ are defined respectively by
\begin{align*}
	& \gph F=\big\{(x, y)\in X\times Y\; \big |\; y\in F(x)\big\}, \; \dom F=\big\{x\in X\; \big |\; F(x)\neq\emptyset\big\}, \\
	&\mbox{\rm and }\rge F=\big\{y\in Y\; \big |\; \exists x\in X\; \mbox{\rm such that }y\in F(x)\big\}.
\end{align*}
We say that $F$ is a {\em convex set-valued mapping} if its graph is a convex set in $X\times Y$. We also say that $F$ is {\em proper} if $\dom F\neq\emptyset$.  We say that $F$ is {\em closed} if $\gph F$ is closed in $X\times Y$.

The {\em epigraphical mapping} $E_f\colon X\tto \R$ of a function $f\colon X\to \oR$, which is defined by
\begin{equation}\label{Emapping}
	E_f(x)=\big\{\lambda \in \R\; \big |\; f(x)\leq \lambda\big\}, \ \; x\in X,
\end{equation} is an important example of a set-valued mapping. It is clear that $\gph E_f=\epi f$ and $\dom E_f=\dom f$. Thus, $f$ is a convex function if and only if $E_f$ is a convex set-valued mapping.

Given a convex set $\Omega\in X$ with $\ox\in \Omega$, define the {\em normal cone} to $\Omega$ at $\ox$ by
\begin{equation}\label{def_normal_cone}
	N(\ox; \Omega)=\big\{x^*\in X^*\; \big |\; \la x^*, x-\ox\ra\leq 0\ \; \mbox{\rm for all }\; x\in \Omega\big\}.
\end{equation}
This concept is used to define the {\em coderivative} of a convex set-valued mapping $F\colon X\tto Y$ at $(\ox, \oy)\in \gph F$ by
\begin{equation}\label{def_coderivative}
	D^*F(\ox, \oy)(y^*)=\big\{x^*\in X^*\; \big |\; (x^*, -y^*)\in N\big((\ox, \oy); \gph F\big)\big\}, \ \; y^*\in Y^*.
\end{equation}
Note that $D^*F(\ox, \oy)$ is a set-valued mapping from $Y^*$ to $X^*$.  The equality
\begin{equation}\label{coder_sub}
	D^*E_f\big(\ox, f(\ox)\big)(1)=\partial f(\ox)
\end{equation} holds for any convex function $f\colon X\to \oR$ and  any $\ox\in X$ with $f(\ox)\in \R$.

Given a function $f\colon X\to \oR$, define the {\em Fenchel conjugate} of $f$ by
\begin{equation*}
	f^*(x^*)=\sup\big\{\la x^*, x\ra-f(x)\; \big |\; x\in X\big\}, \ \; x^*\in X^*.
\end{equation*}
The {\em Fenchel biconjugate} of $f$ is a function $f^{**}\colon X\to \oR$ defined by
\begin{equation*}
	f^{**}(x)=\sup\big\{\la x^*, x\ra-f^*(x^*)\; \big |\; x^*\in X^*\big\}, \ \; x\in X.
\end{equation*}
Among many nice properties of the Fenchel conjugate, we mention the {\em Fenchel-Young inequality} for a proper function $f\colon X\to \oR$:
\begin{equation*}
	\la x^*, x\ra \leq f^*(x^*)+f(x)\ \; \mbox{\rm for all }\; x^*\in X^*,\ x\in X.
\end{equation*}
Furthermore, for $x\in X$ with $f(x)\in\mathbb R$ and $x^*\in X^*$, this inequality holds as an equality if and only if $x^*\in \partial f(x)$.

\section{Fenchel Conjugate of Set-Valued Mappings: Definition and Basic Properties}
\setcounter{equation}{0}

In this section, we define the notion of the Fenchel conjugate for set-valued mappings between topological vector spaces and examine its fundamental properties. We establish the significant connections between this notion and coderivatives of set-valued mappings while generalizing several classical results on the Fenchel conjugate of extended-real-valued functions.

\begin{definition}\label{def_Fenchel_conjugate}
	Let $F\colon X\tto Y$ be a set-valued mapping.  Define $F^*\colon X^*\times Y^*\to \oR$ by
	\begin{equation}\label{Fenchel_conjugate}
		F^*(x^*, y^*)= \sup\big\{\la x^*, x\ra +\la y^*, y\ra \; \big |\; (x, y)\in \gph F\big\}, \ \; (x^*, y^*)\in X^*\times Y^*.
	\end{equation}
	The function $F^*$ is called the  {\em Fenchel conjugate} of $F$.
\end{definition}

By this definition and~\eqref{support function}, the Fenchel conjugate of a set-valued mapping $F\colon X\tto Y$ with nonempty graph is the support function of $\gph F$, i.e.,
\begin{equation}\label{new_label_1}
	F^*(x^*, y^*)=\sigma_{\small \gph F}(x^*, y^*), \ \; (x^*, y^*)\in X^*\times Y^*.
\end{equation}

It is obvious that $F^*$ is a convex function regardless of the convexity of $F$. In addition, if ${\rm gph}\,F\neq\emptyset$, then  by~\cite[Proposition~4.19(a)]{bmn2022} we know that $F^*$ is proper and  lower semicontinuous.

\begin{proposition} Let $F\colon X\tto Y$ be a set-valued mapping.  Then for $x^*\in X^*$ and $y^*\in Y^*$ we have
	\begin{equation*}
		F^*(x^*, 0)=\sigma_{\small\dom F}(x^*)\ \; \mbox{\rm and }\; F^*(0, y^*)=\sigma_{\small\rge F}(y^*).
	\end{equation*}
\end{proposition}
\noindent{\bf Proof.}
	It follows from the definition that
	\begin{align*}
		F^*(x^*, 0)&=\sup\big\{\la x^*, x\ra+\la 0, y\ra\; \big |\; (x, y)\in \gph F\big\}\\
		&=\sup\big\{\la x^*, x\ra\; \big |\; x\in \dom F\big\}=\sigma_{\dom F}(x^*).
	\end{align*}
	The proof for the second formula is similar. $\h$

The proposition below allows us to express the Fenchel conjugate of a function $f\colon X\to\oR$ in terms of the Fenchel conjugate of the epigraphical mapping.
\begin{proposition}\label{fF} Let $f\colon X\to \oR$ be a  function. Then
	\begin{equation*}
		f^*(x^*)=E^*_f(x^*, -1),\ \; x^*\in X^*,
	\end{equation*}
	where $E_f$ is the epigraphical mapping defined in~\eqref{Emapping}.
\end{proposition}
\noindent{\bf Proof.} For any $x^*\in X^*$ we have
	\begin{equation*}
		f^*(x^*)=\sigma_{\epi f}(x^*, -1)=\sigma_{\gph E_f}(x^*, -1)=(E_f)^*(x^*, -1).
	\end{equation*}
	This completes the proof. $\h$

We now discuss the relationships between the Fenchel conjugate and the coderivative of a set-valued mapping.

\begin{theorem}\label{T1} Let $F\colon X\tto Y$ be a  set-valued mapping with $(\ox, \oy)\in \gph F$. For $x^*\in X^*$ and $y^*\in Y^*$ we have the inequality
	\begin{equation}\label{Young1}
		\la x^*, \ox\ra\leq \la y^*, \oy\ra +F^*(x^*, -y^*).
	\end{equation}
	Suppose  in addition that $F$ is convex. Then  $x^*\in D^*F(\ox, \oy)(y^*)$ if and only if
	\begin{equation}\label{Young2}
		\la x^*, \ox\ra= \la y^*, \oy\ra +F^*(x^*, -y^*).
	\end{equation}
\end{theorem}
\noindent{\bf Proof.} Since $(\ox, \oy)\in \gph F$, we have
	\begin{align*}
		\la y^*, \oy\ra+F^*(x^*, -y^*)&=\la y^*, \oy\ra+\sup\{\la x^*, x\ra-\la y^*, y\ra\; |\; (x, y)\in \gph F\}\\
		&\geq \la y^*, \oy\ra+\la x^*, \ox\ra-\la y^*, \oy\ra=\la x^*, \ox\ra.
	\end{align*}
	This implies~\eqref{Young1}.

	Now, assume that $x^*\in D^*F(\ox, \oy)(y^*)$. Then $(x^*, -y^*)\in N((\ox, \oy); \gph F)$ and thus
	\begin{equation*}
		\la x^*, x-\ox\ra-\la y^*, y-\oy\ra\leq 0\ \; \mbox{\rm whenever }\; (x, y)\in \gph F.
	\end{equation*}
	It follows that
	\begin{equation*}
		\la x^*, x\ra-\la y^*, y\ra\leq \la x^*, \ox\ra-\la y^*, \oy\ra\ \; \mbox{\rm for all }(x, y)\in \gph F.
	\end{equation*}
	Taking the supremum with respect to $(x, y)\in \gph F$ gives us
	\begin{equation*}
		F^*(x^*, -y^*)\leq \la x^*, \ox\ra-\la y^*, \oy \ra.
	\end{equation*}
	Combining this with~\eqref{Young1} yields~\eqref{Young2}.
	
	Conversely, if~\eqref{Young2} is valid, then $x^*\in D^*F(\ox, \oy)(y^*)$. The easy verification of this claim is omitted. $\h$

The next example shows that the Fenchel conjugate of a continuous linear mapping can be computed in terms of its adjoint mapping. It follows that the coderivative of the mapping at any point in its graph coincides with the adjoint mapping.

\begin{example}\label{Flinear}{\rm  Let $A\colon X\to Y$ be a continuous linear mapping. The {\em adjoint mapping} $A^*\colon Y^*\to X^*$ of $A$ is defined by $A^*(y^*)=y^*\circ A$  for all $y^*\in Y^*$. In what follows we use $A^*y^*$ to denote  $A^*(y^*)$, and we use $Ax$ to denote $A(x)$ for simplicity. Then we have
		\begin{equation*}
			\la A^*y^*, x\ra=\la y^*, Ax\ra\ \; \mbox{\rm for all }\; x\in X.
		\end{equation*}
		Now, consider the mapping $F(x)=\{Ax\}$ for $x\in X$. Then, by ~\eqref{Fenchel_conjugate}, for $(x^*, y^*)\in X^*\times Y^*$ we have
		\begin{align*}
			F^*(x^*, y^*)&=\sup\{\la x^*, x\ra+\la y^*, Ax\ra\; |\; x\in X\}\\
			&=\sup\{\la x^*+A^*y^*, x\ra\; |\; x\in X\}\\
			&=\begin{cases} 0\; &\mbox{\rm if }A^*y^*=-x^*,\\
				\infty \;& \mbox{\rm otherwise}.
			\end{cases}
		\end{align*}
		Take any $(\ox, \oy)\in \gph F$, i.e., $A\ox=\oy$. By Theorem~\ref{T1}, $x^*\in D^*F(\ox, \oy)(y^*)$ if and only if
		\begin{equation*}
			\la x^*, \ox\ra= \la y^*, \oy\ra +F^*(x^*, -y^*)=\begin{cases} \la y^*, A\ox\ra\; &\mbox{\rm if }A^*y^*=x^*,\\
				\infty \;& \mbox{\rm otherwise}.
			\end{cases}
		\end{equation*}
		This holds if and only if $x^*=A^*y^*$ and thus $D^*F(\ox, \oy)(y^*)=\{A^*y^*\}$.}
\end{example}

Let us show that the well-known results on the Fenchel conjugate for extended-real-valued functions, which have been recalled at the end of Section~\ref{Sect_2}, are consequences of Theorem~\ref{T1} and Proposition~\ref{fF}.

\begin{corollary} Let $f\colon X\to \oR$ be a proper function. Then
	\begin{equation}\label{YF_ineq}
		\la x^*, \ox\ra\leq f^*(x^*)+f(\ox)\ \; \mbox{\rm for all }\; \ox\in X, \;x^* \in X^*.
	\end{equation}
	Assume in addition that $f$ is convex with $\ox\in \dom f$. Then $x^*\in \partial f(\ox)$ if and only if
	\begin{equation*}\label{YF_eq}
		\la x^*, \ox\ra= f^*(x^*)+f(\ox)\ \; \mbox{\rm for all }\; \ox\in X, \; x^* \in X^*.
	\end{equation*}
\end{corollary}
\noindent{\bf Proof.}  Take any  $\ox\in X$ and $x^*\in X^*$. If $f(\bar x)=\infty$, then the inequality in~\eqref{YF_ineq} is obvious. If $f(\bar x)\in \R$, then we put $\bar y=f(\bar x)$ and apply Theorem~\ref{T1} for $F=E_f$ with $y^*:=1$ to get from~\eqref{Young1} the inequality $\langle x^*,\ox\ra\leq \oy+E_f^*(x^*, -1).$ Since $f^*(x^*)=E^*_f(x^*, -1)$  by Proposition~\ref{fF}, this justifies~\eqref{YF_ineq}. Now, suppose that $f$ is convex, $\ox\in \dom f$, and let $\bar y=f(\bar x)$. Then, by Theorem~\ref{T1}, $x^*\in D^*E_f(\ox, \oy)(1)$ if and only if $$\la x^*, \ox\ra= \oy +E_f^*(x^*, -1)$$ or, equivalently,
	$\la x^*, \ox\ra= \oy +f^*(x^*).$ According to~\eqref{coder_sub}, the inclusion $x^*\in D^*E_f(\ox, \oy)(1)$ means that $x^*\in \partial f(\ox)$. So, the second assertion of the corollary is valid. $\h$

\begin{remark}\label{Remark_new1}
	{\rm For any nonempty subset $\Omega\subset X$, one has
		\begin{equation*}\label{support function_property}		\sigma_\Omega(x^*)=\delta_\Omega^*(x^*),\ \; x^*\in X^*.		
		\end{equation*} Indeed, by definitions,
		$$\begin{array}{rcl}
			\sigma_\Omega(x^*) =\sup\big\{\langle x^*,x \rangle\mid x\in \Omega\big\} & = & \sup\big\{\langle x^*,x \rangle-\delta_\Omega(x)\mid x\in X\big\}\\
			& = & \delta_\Omega^*(x^*).
		\end{array}$$}
\end{remark}

The theorem below establishes a relationship between the subdifferential of the Fenchel conjugate of a set-valued mapping and its coderivative.

\begin{theorem} Let $F\colon X\tto Y$ be a proper, closed, convex set-valued mapping. For $(\ox, \oy)\in \gph F$ we have $(\ox, \oy)\in \partial F^*(x^*, y^*)$ if and only if $x^*\in D^*F(\ox, \oy)(-y^*)$.
\end{theorem}
\noindent{\bf Proof.} Fix any $(\ox, \oy)\in \gph F$. Since $\partial F^*(x^*, y^*)=\partial \,\sigma_{\small \gph F}(x^*, y^*)$, thanks to Remark~\ref{Remark_new1} one can rewrite the inclusion  $(\ox, \oy)\in \partial F^*(x^*, y^*)$ as
	\begin{equation}\label{fenchel sym}
		(\ox, \oy)\in \partial h^*(x^*, y^*),
	\end{equation}
	where $h$ is the indicator function associated with $\gph F$, i.e.,
	\begin{equation*}
		h(x, y)=\delta_{\gph F}(x, y), \ \; (x, y)\in X\times Y.
	\end{equation*}
	By~\cite[Theorem~2.4.2(iii)]{zalinescu2002}, we see that~\eqref{fenchel sym}  holds if and only if $$(x^*, y^*)\in \partial h(\ox, \oy)=N\big((\ox, \oy); \gph F\big),$$
	which means that $x^*\in D^*F(\ox, \oy)(-y^*)$. $\h$

\begin{definition}
	The  {\em second Fenchel conjugate} of $F$  is the function $F^{**}\colon X\times Y\to \oR$ given by
	\begin{equation}\label{new_label_2}
		F^{**}(x,y)=(F^*)^*(x, y),\ \; (x, y)\in X\times Y.
	\end{equation}
\end{definition}
\begin{remark}\label{Remark_new2}
	{\rm We have $F^{**}(x, y)=\delta_{\mbox{\rm clco}(\gph F)}(x, y)$, where ${\rm clco}(\Omega)$ denotes the {\em closed convex hull} of a subset $\Omega\subset X\times Y$. Indeed, by~\cite[Proposition~4.19(b)]{bmn2022} one has \begin{equation}\label{equality_new2}\sigma_{\gph F}(x^*,y^*) =\sigma_{\mbox{\rm clco}(\gph F)}(x^*, y^*).\end{equation} Moreover, according to Remark~\ref{Remark_new1}, \begin{equation}\label{equality_new1} \sigma_{\mbox{\rm clco}(\gph F)}(x^*, y^*)= \delta^*_{\mbox{\rm clco}(\gph F)}(x^*, y^*).
		\end{equation} Since $\delta_{\mbox{\rm clco}(\gph F)}(\cdot)$ is a lower semicontinuous and convex  function, by~\cite[Theorem~4.15(c)]{bmn2022} and~\eqref{equality_new1} we assert that $\sigma^*_{\mbox{\rm clco}(\gph F)}(x, y)=  \delta_{\mbox{\rm clco}(\gph F)}(x, y)$. Meanwhile, by~\eqref{new_label_1} we have	$F^*(x^*, y^*)=\sigma_{\small \gph F}(x^*, y^*)$. So, the desired equality $F^{**}(x, y)=\delta_{\mbox{\rm clco}(\gph F)}(x, y)$ follows from~\eqref{new_label_2} and~\eqref{equality_new2}.}
\end{remark}

\begin{proposition} Let $F\colon X\tto Y$ be a set-valued mapping. Then
	\begin{equation*}
		F^{**}(x, y)\leq \delta_{\small\gph F}(x,y) \ \; \mbox{\rm for all }\; (x, y)\in X\times Y.
	\end{equation*}
\end{proposition}
\noindent{\bf Proof.} Since  $\gph F\subset \mbox{\rm clco}(\gph F)$ and $F^{**}(x, y)=\delta_{\mbox{\rm clco}(\gph F)}(x, y)$ by Remark~\ref{Remark_new2}, we obtain the desired result. $\h$

\section{Fenchel Conjugate of  Set-Valued Mappings: Calculus Rules in Finite Dimensions}
\setcounter{equation}{0}

This section is devoted to calculus rules for the Fenchel conjugate of set-valued mappings in finite dimensions. Employing a geometric approach and taking into consideration some relative interior qualification conditions, we establish a sum rule, a chain rule, and an intersection rule for this concept.  Furthermore, we demonstrate that these qualification conditions can be relaxed when polyhedral convexity is present.

Let us begin with a useful representation of the support function associated with a convex polyhedron.

\begin{lemma}\label{lm1} Suppose that $P$ is a nonempty convex polyhedron given by
	\begin{equation*}
		P=\big\{x\in \R^n\; \big |\; \la a_i, x\ra \leq b_i\ \; \mbox{\rm for }\; i=1,\ldots, m\big\}
	\end{equation*} with $a_i\in \R^n$ and $b_i\in \R$ for $i=1,\ldots, m$. Then we have
	$$\sigma_{\small P}(v)= \min\Big\{\sum\limits_{i=1}^m \lambda_ib_i\mid \lambda_1\geq 0, \ldots, \lambda_m\geq 0,\; \sum\limits_{i=1}^m \lambda_ia_i=v \Big\}$$ if $v\in \mbox{\rm pos}\{a_1, \ldots, a_m\}$ and $\sigma_{\small P}(v)=\infty$ if $v\notin \mbox{\rm pos}\{a_1, \ldots, a_m\}$, where $$\mbox{\rm pos}\{a_1, \ldots, a_m\}:=\Big\{\sum\limits_{i=1}^m \lambda_ia_i\mid \lambda_1\geq 0, \ldots, \lambda_m\geq 0\Big\}.$$
\end{lemma}

\noindent{\bf Proof.} By definition, for every $v\in\R^n$, we have
	\begin{equation*}
		\sigma_{\small P}(v)= \sup\big\{\la v,x\ra \mid \la a_i, x\ra \leq b_i\ \; \mbox{\rm for }\; i=1,\ldots, m\big\}.
	\end{equation*}
	Consider the linear programming problem
	\begin{align*}
		&\mbox{\rm maximize }\, \la v, x\ra\\
		&\mbox{\rm subject to }\, \la a_i, x\ra\leq b_i\ \; \mbox{\rm for all }\; i=1, \ldots, m,
	\end{align*}
	which can be re-written as
	\begin{align*}
		&\mbox{\rm maximize }\la v, x\ra\\
		&\mbox{\rm subject to } Ax\leq b\ \; (\mbox{\rm componentwise}),
	\end{align*}
	where $A$ is the matrix whose $i$th row is the row vector $a_i^T$. According to~\cite[Section~2.7]{bss2006} and~\cite[Chapter~5]{Vanderbei2020}, the dual problem of this problem is
	\begin{align*}
		&\mbox{\rm minimize } \la b, \lambda\ra\\
		&\mbox{\rm subject to }A^T\lambda=v,\ \lambda\geq 0.
	\end{align*}
	Since $P$ is nonempty, the primal problem is feasible. If $v\notin \mbox{\rm pos}\{a_1, \ldots, a_m\}$, then the dual problem is infeasible. Hence, by~\cite[Corollary~1, p.~85]{bss2006} we can assert that the primal problem is unbounded, i.e., $\sigma_{\small P}(v)=\infty$. Consider the case where $v\in \mbox{\rm pos}\{a_1, \ldots, a_m\}$. Since both the primal and dual problems are feasible, the strong duality theorem in linear programming (see, e.g.~\cite[Theorem~2.7.3(c)]{bss2006}) assures that both problems have optimal solutions and
	\begin{equation*}
		\sigma_P(v)=\max\big\{\la v,x\ra \mid Ax\leq b\big\} =\min\big\{\la b, \lambda\ra \mid A^T\lambda=v,\ \lambda\geq 0\big\}.
	\end{equation*}
	This implies that
	$$\sigma_{\small P}(v)=\min\Big\{\sum\limits_{i=1}^m \lambda_ib_i\; \big|\; \lambda_1\geq 0, \ldots, \lambda_m\geq 0,\; \sum\limits_{i=1}^m \lambda_ia_i=v \Big\}.$$ This justifies
	the representation of $\sigma_P$ and completes the proof. $\h$

The theorem below allows us to represent the support function associated with the intersection of two convex sets in $\R^n$ in terms of the infimal convolution of the support functions associated with these sets.
\begin{theorem} \label{finite_sigma_intersection_rule} Let $\Omega_1$ and $\Omega_2$ be two nonempty convex sets in $\R^n$. Suppose that one of the following conditions is satisfied:
	\begin{enumerate}
		\item $\Omega_1$ and $\Omega_2$ satisfy the relative interior qualification condition $\ri\Omega_1\cap \ri\Omega_2\neq\emptyset$.
		\item $\Omega_1$ is a polyhedral convex set and
		$\Omega_1\cap \ri\Omega_2\neq\emptyset.$
		\item $\Omega_1$ and $\Omega_2$ are polyhedral convex sets and $\Omega_1\cap \Omega_2\neq\emptyset.$
	\end{enumerate}
	Then $\sigma_{\small \Omega_1\cap \Omega_2}=\sigma_{\small \Omega_1}\s \sigma_{\small \Omega_2}.$
	In addition, for any $v\in \mbox{\rm dom}\big(\sigma_{\small \Omega_1\cap \Omega_2}\big)$, there exist $v_1, v_2\in \R^n$ such that $v=v_1+v_2$ and
	\begin{equation}\label{support rep}
		\sigma_{\small \Omega_1\cap \Omega_2}(v)=\sigma_{\small \Omega_1}(v_1)+\sigma_{\small\Omega_2}(v_2).
	\end{equation}
\end{theorem}
\noindent{\bf Proof.} The proof of the results under (a) and (b) can be found in~\cite[Theorem~4.16]{bmn} and~\cite[Theorem 4.23(c)]{bmn2022}, respectively. Now, we assume that
	\begin{align*}
		&  \Omega_1=\big\{x\in \R^n\; \big|\; \la a_i, x\ra\leq b_i\ \; \mbox{\rm for }\; i=1, \ldots, m\big\},\\
		&\Omega_2=\big\{x\in \R^n\; \big|\; \la a_i, x\ra\leq b_i\ \; \mbox{\rm for }\; i=m+1, \ldots, m+p\big\},
	\end{align*} and the condition $\Omega_1\cap \Omega_2\neq\emptyset$ in~(c) is satisfied. Fix any $v\in \R^n$. It follows from the definition that
	\begin{equation}\label{ineq_inf_conv}
		\sigma_{\small \Omega_1\cap \Omega_2}(v)\leq \big(\sigma_{\small \Omega_1}\s \sigma_{\small \Omega_2}\big)(v).
	\end{equation}
	Thus, if $\sigma_{\small \Omega_1\cap \Omega_2}(v)=\infty$, then the last inequality holds as an equality. Consider the case where $\sigma_{\small \Omega_1\cap \Omega_2}(v)\in \R$, i.e., $v\in \mbox{\rm dom}\big(\sigma_{\small \Omega_1\cap \Omega_2}\big)$. Since $P=\Omega_1\cap \Omega_2$ is a nonempty convex polyhedron, by Lemma~\ref{lm1} we have $v\in \mbox{\rm pos}\{a_1, \ldots, a_{m+p}\}$ and there exist $\lambda_1\geq 0, \ldots, \lambda_m\geq 0$ such that $v=\displaystyle\sum_{i=1}^{m+p}\lambda_i a_i$ and
	\begin{equation*}
		\sigma_{\small \Omega_1\cap \Omega_2}(v)=\sum_{i=1}^{m+p}\lambda_ib_i.
	\end{equation*}
	Setting $v_1=\displaystyle\sum_{i=1}^m \lambda_i a_i$ and $v_2=\displaystyle\sum_{i=m+1}^{m+p} \lambda_i a_i$, we have $$v_1\in \mbox{\rm pos} \{a_1, \ldots, a_m\},\quad v_2\in\mbox{\rm pos}\{a_{m+1}, \ldots,a_{m+p}\},$$ and
	$v=v_1+v_2$. Therefore, using Lemma~\ref{lm1} again gives
	\begin{align*}
		\sigma_{\small \Omega_1\cap \Omega_2}(v)&=\sum_{i=1}^{m}\lambda_ib_i+\sum_{i=m+1}^{m+p}\lambda_ib_i\,  \geq \sigma_{\small \Omega_1}(v_1)+\sigma_{\small \Omega_2}(v_2)\geq \big(\sigma_{\small \Omega_1}\s \sigma_{\small \Omega_2}\big)(v).
	\end{align*}
	Combining this with the inequality~\eqref{ineq_inf_conv}, we can conclude that
	\begin{equation*}
		\sigma_{\small \Omega_1\cap \Omega_2}(v)= \big(\sigma_{\small \Omega_1}\s \sigma_{\small \Omega_2}\big)(v),
	\end{equation*}
	and there exist $v_1, v_2\in \R^n$ with $v=v_1+v_2$ such that~\eqref{support rep} is satisfied. $\h$

The next lemma expresses the Fenchel conjugate of the sum of two set-valued mappings in terms of the support function associated with the intersection of two sets related to the graphs of the set-valued mappings. This lemma plays a crucial role in developing the Fenchel conjugate sum rule for set-valued mappings.
\begin{lemma}\label{lem:conj_sum}
	Let $F_1, F_2\colon \R^n\tto \R^p$ be two set-valued mappings. Define the sets
	\begin{eqnarray}\label{O12}
		\begin{array}{ll}
			&\Omega_1=\big\{(x, y_1, y_2)\in \R^n\times \R^p\times \R^p\; \big|\; y_1\in F_1(x)\big\},\\
			&\Omega_2=\big\{(x, y_1, y_2)\in \R^n\times \R^p\times \R^p\; \big|\; y_2\in F_2(x)\big\}.
		\end{array}
	\end{eqnarray}
	Then we have the equality
	\begin{equation}\label{FC1}
		(F_1+ F_2)^*(u,v)= \sigma_{\Omega_1\cap \Omega_2}(u,v,v)\ \; \mbox{\rm for all }\; (u, v)\in \R^n\times \R^p.
	\end{equation}
\end{lemma}
\noindent{\bf Proof.}  First, consider the case where $\Omega_1\cap \Omega_2\neq\emptyset$. Given any $(u,v)\in \R^n\times \R^p$, we have
	\begin{equation*}
		\begin{aligned}
			\sigma_{\Omega_1\cap \Omega_2}(u,v,v) &= \sup\left\{\la u,x\ra+\la v,y_1\ra+\la v,y_2\ra \mid (x,y_1,y_2)\in\Omega_1\cap\Omega_2\right\}\\
			&= \sup\left\{\la u,x\ra+\la v,y_1+y_2\ra \mid x\in\R^n,\; y_1\in F_1(x),\; y_2\in F_2(x)\right\}\\
			&= \sup\left\{\la u,x\ra+\la v,\ y\ra \mid x\in\R^n,\; y\in (F_1+F_2)(x)\right\}\\
			&=(F_1+F_2)^*(u,v).
		\end{aligned}
	\end{equation*}
	Thus,~\eqref{FC1} is valid. Now, suppose that $\Omega_1\cap \Omega_2=\emptyset$. Then, $\sigma_{\Omega_1\cap \Omega_2}(u,v,v)=-\infty$. Since $\mbox{\rm gph}(F_1+F_2)=\emptyset$, we have $(F_1+ F_2)^*(u,v)=-\infty$. So,~\eqref{FC1} holds. $\h$

The following sum rule is the first main result of this section.

\begin{theorem}\label{Fsumrule}
	Let $F_1, F_2\colon \R^n\tto \R^p$ be two  set-valued mappings and let $(u, v)\in \R^n\times \R^p$. We always have
	\begin{equation}\label{FC2}
		(F_1+F_2)^*(u,v)\leq \inf \left\{F_1^*(u_1,v)+F_2^*(u_2,v)\mid u_1+u_2=u\right\}.
	\end{equation}
	Suppose that one of the following conditions is satisfied:
	\begin{enumerate}
		\item $F_1$ and $F_2$ are convex set-valued mappings such that
		\begin{equation*}\label{QCF}
			\mbox{\rm ri}(\dom F_1)\cap \mbox{\rm ri}(\dom F_2)\neq\emptyset.
		\end{equation*}
		\item $F_1$ is a convex set-valued mapping and $F_2$ is a polyhedral convex set-valued mapping such that
		\begin{equation*}
			\mbox{\rm ri}(\dom F_1)\cap \dom F_2\neq \emptyset.
		\end{equation*}
		\item $F_1$ and $F_2$ are polyhedral convex set-valued mappings such that
		\begin{equation}\label{2poly}
			\dom F_1\cap \dom F_2\neq\emptyset.
		\end{equation}
	\end{enumerate}
	Then we have the equality
	\begin{equation}\label{FC2E}
		(F_1+F_2)^*(u,v)= \inf \left\{F_1^*(u_1,v)+F_2^*(u_2,v)\mid u_1+u_2=u\right\}.
	\end{equation}
	Furthermore, the infimum in~\eqref{FC2E} is attained if $(u, v)\in \mbox{\rm dom}(F_1+F_2)^*$.
\end{theorem}
\noindent{\bf Proof.} Define $\Omega_1$ and $\Omega_2$ as in~\eqref{O12} and fix $(u, v)\in \R^n\times \R^p$. By Lemma~\ref{lem:conj_sum}, we have
	\begin{equation*}
		\begin{aligned}
			(F_1+F_2)^*(u,v)&= 	\sigma_{\Omega_1\cap \Omega_2}(u,v,v) \\
			&=  \sup\left\{\la u,x\ra+\la v,y_1+y_2\ra \mid x\in\R^n,\;  y_1\in F_1(x),\; y_2\in F_2(x)\right\}\\
			&\leq  \sup\left\{\la u_1,x\ra+\la v,y_1\ra \mid x\in\R^n,\; y_1\in F_1(x)\right\}\\&\quad+ \sup\left\{\la u_2,x\ra+\la v,y_2\ra \mid x\in\R^n,\; y_2\in F_2(x)\right\}\\
			&=F_1^*(u_1,v)+F_1^*(u_2,v)
		\end{aligned}
	\end{equation*}
	for any $u_1,u_2\in\R^n$ with $u_1+u_2=u$. Thus, we get~\eqref{FC2}. Now, assume that (a) is satisfied. In this case, $-\infty<(F_1+F_2)^*(u,v)$. If $(F_1+F_2)^*(u,v)=\infty$, then $$\inf \left\{F_1^*(u_1,v)+F_2^*(u_2,v)\mid u_1+u_2=u\right\}=\infty$$ due to~\eqref{FC2}, so~\eqref{FC2} holds as an equality. Next,  consider the situation $(F_1+F_2)^*(u,v)\in \R$, i.e., $(u, v)\in \mbox{\rm dom} (F_1+F_2)^*$. Fix a point $\hat{x}\in \mbox{\rm ri}(\dom F_1)\cap \mbox{\rm ri}(\dom F_2)$. Since $F_1(\hat{x})$ is a nonempty convex set, we can choose $\hat{y}_1\in \ri F_1(\hat{x})$. Similarly, we can choose $\hat{y}_2\in \ri F_2(\hat{x})$. It follows from~\cite[Theorem~2.45]{bmn} that $(\hat x, \hat{y}_1, \hat{y}_2)\in \ri\Omega_1\cap \ri\Omega_2$, and so $\ri\Omega_1\cap \ri\Omega_2\neq\emptyset$. By~\eqref{FC1} we have
	\begin{equation*}
		\sigma_{\Omega_1\cap \Omega_2}(u,v,v)=(F_1+ F_2)^*(u,v)\in \R.
	\end{equation*}
	Using \cite[Theorem~4.16]{bmn} gives us $(u_1, v_1, w_1)$ and $(u_2, v_2, w_2)$ in $\R^n\times \R^p\times \R^p$ such that we have the equalities:
	\begin{eqnarray}\label{FC3}
		\begin{array}{ll}
			&(u, v, v)=(u_1, v_1, w_1)+(u_2, v_2, w_2),\\
			&\sigma_{\Omega_1\cap \Omega_2}(u,v,v)=\sigma_{\Omega_1}(u_1, v_1, w_1)+\sigma_{\Omega_2}(u_2, v_2, w_2).
		\end{array}
	\end{eqnarray}
	If $w_1\neq 0$, then
	\begin{equation*}
		\begin{array}{ll}
			\sigma_{\Omega_1}(u_1, v_1, w_1)&= \sup\left\{\la u_1,x\ra+\la v_1,y\ra+\la w_1,z\ra \mid (x,y)\in\gph F_1,\; z\in\R^p\right\}\\&= \infty.
		\end{array}
	\end{equation*}
	It follows that $\sigma_{\Omega_1\cap \Omega_2}(u,v,v)=\infty$, which is not the case. Thus, we  have $w_1=0$. By a similar argument we deduce that $v_2=0$. Then $u=u_1+u_2$, $v_1=w_2=v$, and by the second equality in~\eqref{FC3} we have
	\begin{align*}
		(F_1+ F_2)^*(u,v)= \sigma_{\Omega_1\cap \Omega_2}(u,v,v)&=\sigma_{\gph F_1}(u_1, v)+\sigma_{\gph F_2}(u_2, v)\\
		&=F^*_1(u_1, v)+F^*_2(u_2, v)\\
		&\geq \inf \left\{F_1^*(u_1,v)+F_2^*(u_2,v)\mid u_1+u_2=u\right\}.
	\end{align*}
	Combining this with~\eqref{FC2} gives~\eqref{FC2E}. Our arguments also justify that the infimum in~\eqref{FC2E} is attained if $(u, v)\in \mbox{\rm dom}(F_1+F_2)^*$.
	
	If (b) is satisfied, then we can find $\hat{x}\in \mbox{\rm ri}(\dom F_1)\cap \dom F_2$. Since $F_1(\hat{x})$ is a nonempty convex set in $\R^n$, there exists $\hat{y}_1\in \mbox{\rm ri}\, F_1(\hat{x}_1)$. By \cite[Theorem~2.45]{bmn},~$(\hat{x}, \hat{y}_1)\in \mbox{\rm ri}(\gph F_1)$. Then, choosing $\hat{y}_2\in F_2(\hat{x})$, we get
	$(\hat{x}, \hat{y}_1, \hat{y}_2)\in \mbox{\rm ri}(\Omega_1)\cap \Omega_2.$ Since~$\Omega_2$ is a polyhedral convex set, repeating the proof of the theorem under case~(a) with the use of Theorem~\ref{finite_sigma_intersection_rule}(b) we get the desired assertions.
	
	Finally, assume that (c) is satisfied. Then both $\Omega_1$ and $\Omega_2$ are polyhedral convex sets. By~\eqref{2poly} we have  $\Omega_1\cap \Omega_2\neq\emptyset$. Therefore, we can apply Theorem~\ref{finite_sigma_intersection_rule}(c) to complete the proof in this case. $\h$

The lemma below expresses the Fenchel conjugate of the composition of two set-valued mappings in term of the intersection of two sets whose definitions include the graphs of the set-valued mappings.
\begin{lemma}
	\label{lem:chainrule}
	Let $F\colon\R^n\tto \R^p$ and $G:\R^p\tto \R^q$ be arbitrary set-valued mappings, and let
	\begin{eqnarray}\label{omegas}
		\Omega_1=(\gph F)\times \R^q,\ \;\Omega_2=\R^n\times(\gph G).
	\end{eqnarray}
	Then $(G\circ F)^*(u,w) =\sigma_{\Omega_1\cap\Omega_2}(u,0,w)$ for every $(u,w)\in \R^n\times\R^q$.
\end{lemma}
\noindent{\bf Proof.} Take any $(x,y,z)\in \Omega_1\cap \Omega_2$. Then, $z\in (G\circ F) (x)$. For any $(u,w)\in\R^n\times\R^q $, we have
	\begin{equation*}
		\begin{aligned}
			\la u, x\ra + \la 0,y\ra + \la w,z\ra &= \la u, x\ra + \la w,z\ra \\
			&\leq \sup\left\{\la u,x'\ra +\la w,z'\ra \mid (x',z')\in \mbox{\rm gph}(G\circ F)\right\}\\
			&= (G\circ F)^*(u,w).
		\end{aligned}
	\end{equation*}
	We thus can deduce that
	\begin{equation}\label{chain1}
		\sigma_{\Omega_1\cap\Omega_2}(u,0,w) \leq  (G\circ F)^*(u,w).
	\end{equation}
	Let us now take any $x\in\R^n$ and $z\in (G\circ F)(x)$. Then there exists $y\in \R^p$ such that $z\in G(y)$ and $y\in F(x)$. Hence, $(x,y,z)\in \Omega_1\cap\Omega_2$ and we get
	\begin{align*}
		\la u,x\ra +\la w,z\ra &= 	\la u, x\ra + \la 0,y\ra + \la w,z\ra\\
		&\leq \sup\left\{\la u,x'\ra +\la 0,y'\ra+\la w,z'\ra \mid (x',y',z')\in\Omega_1\cap \Omega_2\right\}\\
		&= 	\sigma_{\Omega_1\cap\Omega_2}(u,0,w).
	\end{align*}
	Since the latter inequality holds for all $(x,z)\in \mbox{\rm gph}(G\circ F)$, we have
	\begin{equation*}
		(G\circ F)^*(u,w) = \sup \left\{\la u,x\ra +\la w,z\ra \mid z\in (G\circ F)(x) \right\} \leq \sigma_{\Omega_1\cap\Omega_2}(u,0,w).
	\end{equation*}
	Together with~\eqref{chain1}, this yields the desired result.
$\h$

The next Fenchel conjugate chain rule for set-valued mappings is another main result of this section.

\begin{theorem}\label{Fchainrulef}
	Let $F\colon \R^n\tto \R^p$ and $G\colon \R^p\tto \R^q$ be two set-valued mappings. Then for any $(u,w)\in \R^n\times\R^q$ we have
	\begin{equation}\label{Fchainruleff}
		(G\circ F)^*(u,w) \leq \inf\{F^*(u,v)+G^*(-v,w)\mid v\in\R^p\}.
	\end{equation}
	The equality
	\begin{equation}\label{FchainruleffE}
		(G\circ F)^*(u,w) = \inf\{F^*(u,v)+G^*(-v,w)\mid v\in\R^p\}
	\end{equation} holds under any one of the following conditions:
	\begin{enumerate}
		\item $F$ and $G$ are convex set-valued mappings such that
		\begin{equation}\label{QCC}
			\mbox{\rm ri}(\rge F)\cap \mbox{\rm ri}(\dom G)\neq\emptyset.
		\end{equation}
		\item $F$ is a polyhedral convex set-valued mapping and $G$ is a convex set-valued mapping such that $
		\rge F\cap \mbox{\rm ri}(\dom G)\neq\emptyset$ (or $F$ is a convex set-valued mapping and $G$ is a polyhedral convex set-valued mapping such that $\mbox{\rm ri}(\rge F)\cap \dom G\neq\emptyset$).
		\item Both $F$ and $G$ are polyhedral convex set-valued mappings such that
		\begin{equation*}
			\rge F\cap \dom G\neq\emptyset.
		\end{equation*}
	\end{enumerate}
	Furthermore, in that case, the infimum in~\eqref{FchainruleffE} is attained if $(G\circ F)^*(u,w)\in \R$.
\end{theorem}
\noindent{\bf Proof.} Take any $x\in\R^n$ and $z\in (G\circ F)(x)$. Then there exists $y\in \R^p$ such that $z\in G(y)$ and $y\in F(x)$. For any $(u,w)\in\R^n\times \R^q$ and $v\in\R^p$, we have
	\begin{align*}
		\la u,x \ra +\la w,z\ra &=  	\la u,x \ra +\la v,y\ra + \la -v,y\ra + \la w,z\ra\\
		&\leq \sup \left\{	\la u,x' \ra +\la v,y'\ra \mid (x',y')\in \gph F\right\} \\&\quad+ \sup \left\{	\la -v,y' \ra +\la w,z'\ra \mid (y',z')\in \gph G\right\}\\
		&= F^*(u,v)+G^*(-v,w).
	\end{align*}	
	Therefore, for any $(u,w)\in\R^n\times\R^q$,
	\begin{equation*}\label{chain2}
		\begin{aligned}
			(G\circ F)^*(u,w) &=\sup \left\{\la u,x\ra +\la w,z\ra \mid z\in (G\circ F)(x) \right\} \\
			&\leq \inf\left\{F^*(u,v)+G^*(-v,w)\mid v\in\R^p\right\},
		\end{aligned}
	\end{equation*} so~\eqref{Fchainruleff} holds.
	To prove the reverse inequality, let us define $\Omega_1$ and $\Omega_2$ as in~\eqref{omegas}. Clearly, the convexity of $F$ and $G$ implies that $\Omega_1$ and $\Omega_2$ are convex.
	
	Since~\eqref{QCC} yields $(G\circ F)^*(u,w) >-\infty$ for all $(u,w)\in\R^n\times\R^q,$ it follows from Lemma~\ref{lem:chainrule} that $\sigma_{\Omega_1\cap\Omega_2}(u,0,w) >-\infty$ for all $(u,w)\in\R^n\times\R^q$. If for some $(u,w)\in\R^n\times\R^q$ one has $\sigma_{\Omega_1\cap\Omega_2}(u,0,w) =\infty$, then  $(G\circ F)^*(u,w) =\infty$. So, the inequality~\eqref{Fchainruleff} implies that
	\begin{equation*}
		\infty = (G\circ F)^*(u,w) = \inf\left\{F^*(u,v)+G^*(-v,w)\mid v\in\R^p\right\}.
	\end{equation*}
	It remains to consider all $(u,w)\in \R^n\times\R^q$ such that $(u,0,w)\in\mbox{\rm dom} (\sigma_{\Omega_1\cap\Omega_2})$. From~\cite[Theorem~2.45]{bmn} and the condition~\eqref{QCC}, it follows that
	\begin{equation}\label{ri_intersection}
		\ri \Omega_1\cap\ri\Omega_2\neq\emptyset.
	\end{equation} Indeed, we have $\dom (F^{-1})=\rge (F)$. Since $T(x, y):=(y, x)$ for $(x,y)\in X\times Y$ is a linear homeomorphism between $X\times Y$ and $Y\times X$, one has
	$\ri T(A)=T(\ri A)$ for any convex subset $A\subset X\times Y$. Note that $T(\gph F)=\gph F^{-1}$. By~\eqref{QCC}, we can find a point $y_0 \in \ri(\rge F) \cap \ri(\dom G)$. Pick any $x_0\in \ri F^{-1}(y_0)$ and $z_0\in\ri G(y_0)$ (recall that the relative interior of a nonempty convex subset of a finite-dimensional Euclidean space has a nonempty relative interior). Then $(y_0, z_0)\in \ri(\gph (G))$ by~\cite[Theorem~2.45]{bmn}, so $(x_0, y_0, z_0)\in \ri\Omega_2$. Furthermore, by~\cite[Theorem~2.45]{bmn} we have $$(y_0, x_0)\in \ri(\gph F^{-1})=\ri T(\gph F)=T(\ri(\gph F)).$$ It follows that $(y_0, x_0)=T(x_1, y_1)=(y_1, x_1)$ for some element $(x_1, y_1)\in \ri(\gph F)$. Thus $(x_0, y_0)\in \ri(\gph F)$, so
	$(x_0, y_0, z_0)\in \ri(\Omega_1)$. Then the inclusion $(x_0, y_0, z_0)\in \ri(\Omega_1)\cap \ri(\Omega_2)$ is valid. So,~\eqref{ri_intersection} holds.

	Thanks to~\eqref{ri_intersection} and the property $(u,0,w)\in\mbox{\rm dom} (\sigma_{\Omega_1\cap\Omega_2})$, applying~\cite[Theorem~4.16]{bmn} for $\Omega_1$ and $\Omega_2$,   we can find $(u_1,v,w_1),\ (u_2,-v,w_2)\in \R^n\times\R^p\times\R^q$ such that
	\begin{equation}\label{u_w}
		u= u_1+u_2,\quad w = w_1+w_2,
	\end{equation}
	and
	\begin{equation}
		\label{chain3}
		\sigma_{\Omega_1\cap \Omega_2}(u,0,w) = \sigma_{\Omega_1}(u_1,v,w_1) +\sigma_{\Omega_2}(u_2,-v,w_2).
	\end{equation}
	If $w_1\neq 0$, then $\sigma_{\Omega_1}(u_1, v, w_1)= \infty$, which together with~\eqref{chain3} contradicts the fact that $\sigma_{\Omega_1\cap \Omega_2}(u,0,w) \in \R$. Similarly, we get $u_2=0$. Hence, from~\eqref{u_w} it follows that $u_1=u$ and $w_2=w$. Thus, by Lemma~\ref{lem:chainrule} and the equality~\eqref{chain3} one has
	\begin{equation*}
		\begin{aligned}
			(G\circ F)^*(u,w)= \sigma_{\Omega_1\cap \Omega_2}(u,0,w) &\,= \sigma_{\Omega_1}(u,v,0) +\sigma_{\Omega_2}(0,-v,w) \\&= \sigma_{\gph F}(u,v)+\sigma_{\gph G}(-v,w)\\
			&=F^*(u,v)+G^*(-v,w)\\
			&\geq  \inf\left\{F^*(u,\tilde v)+G^*(-\tilde v,w)\mid \tilde v\in\R^p\right\}.
		\end{aligned}
	\end{equation*}
	Combining this with~\eqref{Fchainruleff} yields
	\begin{equation*}
		(G\circ F)^*(u,w) = \inf\left\{F^*(u,v)+G^*(-v,w)\mid v\in\R^p\right\}
	\end{equation*}
	for all $(u,w)\in \R^n\times\R^q$ such that $(u,0,w)\in\mbox{\rm dom} (\sigma_{\Omega_1\cap\Omega_2})$. The proof of the theorem under~(b) or~(c)  is based on the above arguments and the assertions of Theorem~\ref{finite_sigma_intersection_rule} under the conditions~(b) or~(c) there.	$\h$

We conclude this section with an intersection rule for the Fenchel conjugate of  set-valued mappings.
\begin{theorem}
	Let $F_1, F_2\colon \R^n\tto \R^p$ be two convex set-valued mappings. Then we have the inequality
	\begin{equation*}
		(F_1\cap F_2)^*(u, v)\leq \big(F_1^*\square F_2^*\big)(u, v)\ \; \mbox{\rm for all }\; (u, v)\in \R^n\times \R^p.
	\end{equation*}
	The equality
	\begin{equation}\label{FintersectionE}
		(F_1\cap F_2)^*(u, v)= \big(F_1^*\square F_2^*\big)(u, v)\ \; \mbox{\rm for all }\; (u, v)\in \R^n\times \R^p
	\end{equation}
	holds under one of the following conditions:
	\begin{enumerate}
		\item $F_1$ and $F_2$ are convex set-valued mappings such that
		$\mbox{\rm ri}(\gph F_1)\cap \mbox{\rm ri}(\gph F_2)\neq\emptyset$.
		\item $F_1$ is a convex set-valued mapping and $F_2$ is a polyhedral convex set-valued mapping such that
		$\gph F_1\cap \mbox{\rm ri}(\gph F_2)\neq\emptyset$.
		\item Both $F_1$ and $F_2$ are polyhedral convex set-valued mappings such that
		\begin{equation*}
			\gph F_1\cap \gph F_2\neq\emptyset.
		\end{equation*}
	\end{enumerate}
	Furthermore, if $(F_1\cap F_2)^*(u, v)\in \R$, then the infimum in the infimal convolution
	\begin{equation*}
		\big(F_1^*\square F_2^*\big)(u, v)=\inf\big\{F_1^*(u_1,v_1)+F_2^*(u_2,v_2)\; \big|\; (u_1,v_1)+(u_2,v_2)=(u,v)\big\}
	\end{equation*}
	on the right-hand-side of~\eqref{FintersectionE} is attained.
\end{theorem}
\noindent{\bf Proof.} It follows from the definition that
	\begin{equation*}
		(F_1\cap F_2)^*(u, v)=\sigma_{\mbox{\rm gph}(F_1\cap F_2)}(u,v) = \sigma_{(\gph F_1)\cap (\gph F_2)}(u, v),
	\end{equation*}
	for any $(u,v)\in\R^n\times\R^p$.Then we only need to apply Theorem~\ref{finite_sigma_intersection_rule} to complete the proof.
$\h$

\section{Fenchel Conjugate of Set-Valued Mappings: Calculus Rules in Infinite Dimensions}
\setcounter{equation}{0}

This section establishes some useful calculus rules for the Fenchel conjugate of convex set-valued mappings  in locally convex topological vector spaces and in Banach spaces.  The rules require qualifications conditions that rely on the concepts of interiors and generalized relative interiors. Similar to the finite-dimensional case, the qualification conditions can be relaxed when polyhedral convexity is present.

The \textit{quasi-relative interior} and the  \textit{strong quasi-relative interior} of a subset $\Omega$ of a locally convex topological vector space $X$ are defined (see~\cite{bmn2022} and the references therein) respectively by
\begin{align*}
	&\mbox{\rm qri}\, \Omega=\big\{a\in \Omega\; |\; \overline{\mbox{\rm cone}(\Omega-a)}\; \mbox{\rm is a linear subspace of }X\big\},\\
	&\mbox{\rm sqri}\, \Omega=\big\{a\in \Omega\; |\; \mbox{\rm cone}(\Omega-a)\; \mbox{\rm is a closed  linear subspace of }X\big\}.
\end{align*} It is clear that $\mbox{\rm sqri}\, \Omega\subset \mbox{\rm qri}\, \Omega$. The theorem below provides us with a support function intersection rule in locally convex topological vector spaces under different settings. While a significant portion of the results can be found in~\cite{bmn2022}, we  include additional details of the proof for the reader's convenience.
\begin{theorem}\label{Flemma1} Let $\Omega_1$ and $\Omega_2$ be two nonempty convex sets in $X$. Suppose that one of the following conditions is satisfied:
	\begin{enumerate}
		\item $\Omega_1\cap \mbox{\rm int}(\Omega_2)\neq\emptyset$.
		\item $\Omega_1$ is a polyhedral convex set and $\Omega_1\cap (\qri \Omega_2)\neq\emptyset$.
		\item $\Omega_1$ and $\Omega_2$ are polyhedral convex sets and $\Omega_1\cap \Omega_2\neq\emptyset$.
		\item $\Omega_1$ and $\Omega_2$ are closed sets, $X$ is a Banach space, and $0\in\mbox{\rm sqri}\,(\Omega_1-\Omega_2)$.
	\end{enumerate}
	Then we have $\sigma_{\Omega_1\cap \Omega_2}=\sigma_{\Omega_1}\s \sigma_{\Omega_2}$.
	Furthermore, if $\big(\sigma_{\Omega_1\cap \Omega_2}\big)(x^*)\in \R$, there exist $x^*_1, x^*_2\in X^*$ with $x^*=x^*_1+x^*_2$ and
	\begin{equation*}
		\big(\sigma_{\Omega_1\cap \Omega_2}\big)(x^*)=\sigma_{\Omega_1}(x^*_1)+\sigma_{\Omega_2}(x^*_2).
	\end{equation*}
\end{theorem}

\noindent{\bf Proof.} Fix $x^*\in X^*$ and pick $x_1^*$ and $x_2^*$ with $x^*=x_1^*+x_2^*$. Then for every $x\in\Omega_1\cap\Omega_2$, we have
	\begin{equation*}
		\la x^*,x\ra = 	\la x_1^*,x\ra+ \la x_2^*,x\ra \leq \sigma_{\Omega_1}(x_1^*)+ \sigma_{\Omega_2}(x_2^*).
	\end{equation*}
	Then, taking the infimum on the second sum of the above expression with respect to $x_1^*$ and $x_2^*$ yields
	\begin{equation*}
		\la x^*,x\ra \leq (\sigma_{\Omega_1}\s \sigma_{\Omega_2})(x^*), \ \,\text{for all}\; x\in\Omega_1\cap\Omega_2.
	\end{equation*}
	Since $x^*$ can be chosen arbitrarily, $\sigma_{\Omega_1\cap \Omega_2}(x^*)\leq\big(\sigma_{\Omega_1}\s \sigma_{\Omega_2}\big)(x^*)$. It remains to prove that, under each one of the conditions (a)--(d), the inequality
	\begin{equation}\label{support_intersection_geq}
		\big(\sigma_{\Omega_1\cap \Omega_2}\big)(x^*)\geq \big(\sigma_{\Omega_1}\s \sigma_{\Omega_2}\big)(x^*)
	\end{equation}
	holds for every $x^*\in X^*$.\\[1ex]
	(a)  If $x^*\notin \mbox{\rm dom}(\sigma_{\Omega_1\cap \Omega_2})$, then~\eqref{support_intersection_geq} is obvious. Fix any  $x^*\in \mbox{\rm dom}(\sigma_{\Omega_1\cap \Omega_2})$ and denote
	\begin{equation}\label{alpha}
		\alpha =\sigma_{\Omega_1\cap \Omega_2}(x^*).
	\end{equation}
	Let
	\begin{eqnarray}\label{thetas}
		\begin{array}{ll}
			&\Theta_1= \left\{(x,\lambda)\in X\times\R\mid x\in\Omega_1,\; \lambda \leq \la x^*,x\ra -\alpha\right\},\\
			&\Theta_2=\Omega_2\times [0,\infty).
		\end{array}
	\end{eqnarray}
	Since $\mbox{\rm int}(\Omega_2)\neq\emptyset$, we have $\mbox{\rm int}(\Theta_2) \neq \emptyset$. In addition, from~\eqref{alpha} and~\eqref{thetas} it follows that $\Theta_1\cap \mbox{\rm int}(\Theta_2)=\emptyset$. Therefore, by~\cite[Corollary~2.59]{bmn2022} we can find a nonzero vector $(y^*,\beta)\in X^*\times\R$ such that
	\begin{equation}
		\label{separate_thetas}
		\la y^*,y\ra +\beta\lambda_2\leq \la y^*,x\ra+\beta\lambda_1,
	\end{equation}
	for any $(x,\lambda_1)\in\Theta_1$ and $(y,\lambda_2)\in\Theta_2$. Then, from the definitions of $\Theta_1$ and $\Theta_2$, it follows that $\beta \leq 0$. If $\beta = 0$, then we can deduce from~\eqref{separate_thetas} that
	\begin{equation}\label{separation_0}
		\la y^*,y-x\ra \leq 0,\ \, \text{for all}\; x\in\Omega_1,\; y\in\Omega_2.
	\end{equation}
	By our assumption, there exists $\bar x\in \Omega_1\cap \mbox{\rm int}(\Omega_2)$. So, we have
	$$0\in \mbox{\rm int}(\Omega_1)-\bar x\subset \mbox{\rm int}(\Omega_1-\Omega_2),$$ where the second inclusion holds because $\mbox{\rm int}(\Omega_1)-\bar x$ is an open subset of $\Omega_1-\Omega_2$. Hence,~\eqref{separation_0} implies that $y^*=0$, which is a contradiction as $(y^*,\beta)$ is nonzero. So, $\beta <0$. For any $(x,y)\in\Omega_1\times\Omega_2$, we have $(x,\la x^*,x\ra -\alpha)\in \Theta_1$ and $(y,0)\in\Theta_2$, then from~\eqref{separate_thetas} it follows that
	\begin{equation*}
		\la y^*,y\ra \leq  \la y^*,x\ra + \beta (\la x^*,x\ra -\alpha)\ \; \text{for all}\ \; x\in\Omega_1,\; y\in\Omega_2,
	\end{equation*}
	or equivalently,
	\begin{equation*}
		\beta\alpha\leq \la y^*+\beta x^*, x\ra +\la -y^*,y\ra\ \; \text{for all}\ \; x\in\Omega_1,\; y\in\Omega_2.
	\end{equation*}
	Dividing both sides of the last inequality by $\beta$ and noting that $\beta <0$, we have
	\begin{equation*}
		\alpha \geq  \left\la \frac{y^*}{\beta}+x^*, x\right\ra +\left\la \frac{-y^*}{\beta},y\right\ra\ \; \text{for all}\ \; x\in\Omega_1, y\in\Omega_2.
	\end{equation*}
	By putting $x_1^*= y^*/\beta+x^*$ and $x_2^*=-y^*/\beta$, one obtains~\eqref{support_intersection_geq}.\\[1ex]
	(b)  Again, define $\alpha$ as in~\eqref{alpha} and $\Theta_1$ and $\Theta_2$ as in~\eqref{thetas}. Since $\Omega_1$ is a polyhedral convex set, $\Theta_1$ is also a polyhedral convex set. Clearly, $\Theta_2$ is a convex set, with $\qri \Theta_2=\qri \Omega_2\times (0,\infty)$. Since $\Omega_1\cap (\qri \Omega_2)\neq\emptyset$, we have $\qri \Theta_2\neq \emptyset$. By the definition of $\alpha$, it is easy to verify that $\Theta_1\cap (\qri \Theta_2)=\emptyset$. Then, applying~\cite[Theorem~3.86]{bmn2022} we have that $\Theta_1$ and $\Theta_2$ can be separated by a closed hyperplane that does not contain $\Theta_2$. That is, there exist a nonzero vector $(y^*,\beta)\in X^*\times\R$ and a number $\gamma\in \R$ such that
	\begin{equation}
		\label{separate_thetas_2}
		\la y^*,y\ra +\beta\lambda_2\leq \gamma \leq \la y^*,x\ra+\beta\lambda_1,
	\end{equation}
	for any $(x,\lambda_1)\in\Theta_1$ and $(y,\lambda_2)\in\Theta_2$.
	In addition, since the separating hyperplane
	$$H=\{(x, \lambda)\in X\times \R\; |\; \la y^*, x\ra+\beta \lambda=\gamma\}$$
	does not contain $\Theta_2$, there exists $(\hat{y}, \hat{\lambda}_2)\in \Theta_2$ such that
	\begin{equation}\label{theta2_notin_H}
		\la y^*,\hat{y}\ra +\beta\hat{\lambda}_2< \gamma.
	\end{equation}
	As above, the definitions of $\Theta_1$ and $\Theta_2$ give that $\beta \leq 0$. If $\beta = 0$, then~\eqref{separate_thetas_2} reduces to
	\begin{equation}\label{separate_lambdas}
		\la y^*,y\ra \leq \gamma \leq \la y^*,x\ra\; \text{ for any }\; x\in\Omega_1,\ y\in\Omega_2
	\end{equation}
	and~\eqref{theta2_notin_H} reduces to
	\begin{equation}\label{Omega2_notin_H}
		\la y^*, \hat y\ra < \gamma\; \text{ with }\;  \hat y\in\Omega_2.
	\end{equation}
	Now~\eqref{separate_lambdas} shows that the polyhedral set $\Omega_1$ and the convex set $\Omega_2$ can be separated by the closed hyperplane $H_1=\{x\in X\; |\; \la y^*, x\ra=\gamma\}$, and~\eqref{Omega2_notin_H} shows both that the separation is proper and that $\Omega_2$ is not contained in $H_1$. We again apply~\cite[Theorem~3.86]{bmn2022} to get $\Omega_1\cap(\qri\Omega_2)=\emptyset$  which is a contradiction, thus $\beta<0$. We now conclude as in part~(a), obtaining~\eqref{support_intersection_geq}.\\[1ex]
	(c)  Since $\Omega_1$ and $\Omega_2$ are polyhedral convex sets, there exist $x^*_i\in X^*$ and $b_i\in\R$ for $i=1,\dots,m+p$ and $m,p\in\mathbb{N}$ such that $\Omega_1$ and $\Omega_2$ admit the following representations
	\begin{eqnarray*}
		\begin{array}{ll}
			&\Omega_1=\{x\in X\mid x^*_i(x)\leq b_i,\; i=1,\dots,m\},\\
			&\Omega_2=\{x\in X\mid x^*_i(x)\leq b_i,\; i=m+1,\dots,m+p\}.
		\end{array}
	\end{eqnarray*}
	Let $L= \bigcap_{i=1}^{m+p} \mbox{\rm ker}x^*_i$. For each $i\in\{1,\dots,m+p\}$, define the functions $\hat{x}^*_i\colon X/L\to\R$ by
	\begin{equation*}
		\hat{x}^*_i([x])= x^*_i(x),\ \; \text{where}\ \; [x]=x+L\in X/L.
	\end{equation*}
	We see that the function $\hat{x}^*_i$ is well-defined and belongs to $(X/L)^*$ for every $i=1,\dots,m+p$. Consider the quotient map $\pi\colon X\to X/L$ defined by
	$\pi(x)=x+L$ for $x\in X$. Then we have
	\begin{equation*}
		\pi(\Omega_1)= \big\{[x]\in X/L\mid \hat{x}^*_i([x])\leq b_i,\; i=1,\dots,m\big\}
	\end{equation*}
	and
	\begin{equation*}
		\pi(\Omega_2)= \big\{[x]\in X/L\mid \hat{x}^*_i([x])\leq b_i,\; i=m+1,\dots,m+p\big\}.
	\end{equation*}
	So, $\pi(\Omega_1)$ and $\pi(\Omega_2)$ are polyhedral convex sets in $X/L$. Since $\Omega_1\cap\Omega_2\neq \emptyset$, we see that $\pi(\Omega_1)\cap \pi(\Omega_2)\neq \emptyset$. It follows from~\cite[Proposition~1.119]{bmn2022} that the quotient space $X/L$ is finite-dimensional. Let us denote by $\pi^*\colon (X/L)^*\to X^*$ the adjoint operator of $\pi$. By~\cite[Proposition~1.118]{bmn2022}, $\pi^*$ is onto the space $L^\perp$. Note that we only need to prove the inequality~\eqref{support_intersection_geq} for the case where $\alpha=(\sigma_{\Omega_1\cap \Omega_2})(x^*)\in \R$. In this situation, $x^*\in L^\perp$. Indeed, suppose on the contrary that $x^*\notin L^\perp$. Then we find $x_0\in L$ such that $\la x^*, x_0\ra>0$. Choose $\ox\in \Omega_1\cap \Omega_2$ and set $x_t=\ox+tx_0$ for $t>0$. Since $x_t\in \Omega_1\cap \Omega_2$ for all $t>0$, we get
	\begin{equation*}
		\sigma_{\Omega_1\cap \Omega_2}(x^*)\geq \sup_{t>0}\,\la x^*, \ox+tx_0\ra=\infty,
	\end{equation*}
	which is a contradiction.
	
	Choose $f\in (X/L)^*$ such that $\pi^* f=x^*$ and get
	$\la x^*, x\ra=\la \pi^* f, x\ra=\la f, [x]\ra$ for all $x\in X$. Then
	\begin{align*}
		\sigma_{\Omega_1\cap \Omega_2}(x^*)&= \sup\{\la x^*,x\ra\mid x\in\Omega_1\cap\Omega_2\}\\
		&= \sup\{\la f,[x]\ra \mid [x]\in \pi(\Omega_1)\cap \pi(\Omega_2)\} \\
		&= \sigma_{\pi(\Omega_1)\cap \pi(\Omega_2)}(f).
	\end{align*}
	By Theorem~\ref{finite_sigma_intersection_rule}(c), there exist $f_1, f_2\in \ (X/L)^*$ such that $f=f_1+f_2$ and
	\begin{equation*}
		\sigma_{\pi(\Omega_1)\cap \pi(\Omega_2)}(f)=\sigma_{\pi(\Omega_1)}(f_1)+\sigma_{\pi(\Omega_1)}(f_2).
	\end{equation*}
	Denote $x^*_1=\pi^*f_1\in X^*$ and $x^*_2=\pi^*f_2\in X^*$. Then $x^*=x^*_1+x^*_2$ and
	\begin{align*}
		\sigma_{\pi(\Omega_1)\cap \pi(\Omega_2)}(f)&=\sigma_{\pi(\Omega_1)}(f_1)+\sigma_{\pi(\Omega_1)}(f_2)\\
		&= \sup_{x_1\in \Omega_1}\la f_1, [x_1]\ra +\sup_{x_2\in \Omega_2}\la f_2, [x_2]\ra\\
		& = \sup_{x_2\in \Omega_1}\la x^*_1, x_1\ra+\sup_{x_2\in \Omega_2}\la x^*_2, x_2\ra\ \\
		&=\sigma_{\Omega_1}(x^*_1)+\sigma_{\Omega_2}(x^*_2)
		\geq (\sigma_{\Omega_1}\s \sigma_{\Omega_2})(x^*).
	\end{align*}
	(d) The conclusion under (d) is well-known; see, e.g., \cite{AB}. $\h$

Based on Theorem~\ref{Flemma1}, we obtain below the Fenchel conjugate sum rule for set-valued mappings in infinite dimensions.
\begin{theorem}\label{Fsumrulei} Let $F_1, F_2\colon X\tto Y$ be set-valued mappings between locally convex topological vector spaces and let $(x^*, y^*)\in X^*\times Y^*$. We always have
	\begin{equation}\label{Fsumruleif}
		(F_1+F_2)^*(x^*,y^*)\leq \inf \left\{F_1^*(x^*_1,y^*)+F_2^*(x^*_2,y^*)\mid x^*_1+x^*_2=x^*\right\}.
	\end{equation}
	Suppose that one of the following conditions is satisfied:
	\begin{enumerate}
		\item  $F_1$ and $F_2$ are convex set-valued mappings such that $\mbox{\rm int}(\gph F_1)\neq\emptyset$ and
		\begin{equation*}\label{QCF_2}
			\mbox{\rm int}(\dom F_1)\cap \dom F_2\neq\emptyset.
		\end{equation*}
		\item $F_1$ is a  convex set-valued mapping and $F_2$ is a polyhedral convex set-valued mapping such that
		there exist $\hat{x}\in\mbox{\rm qri}(\dom F_1)\cap \dom F_2$ and $\hat{y}_1\in Y$ with $(\hat{x}, \hat{y}_1)\in \mbox{\rm qri}(\gph F_1)$.
		\item $F_1$ and $F_2$ are polyhedral convex set-valued mappings and
		$\dom F_1\cap \dom F_2\neq \emptyset$.
		\item $X$ and $Y$ are Banach spaces, $F_1$ and $F_2$ have closed convex graphs, and \begin{equation}\label{regu_d} 0\in \mbox{\rm sqri}(\dom F_1-\dom F_2).\end{equation}
	\end{enumerate}
	Then we have the equality
	\begin{equation}\label{FsumruleifE}
		(F_1+F_2)^*(x^*,y^*)= \inf \left\{F_1^*(x^*_1,y^*)+F_2^*(x^*_2,y^*)\mid x^*_1+x^*_2=x^*\right\}.
	\end{equation}
	Furthermore, the infimum is attained in~\eqref{FsumruleifE} if $(x^*, y^*)\in \mbox{\rm dom}(F_1+F_2)^*$.
\end{theorem}
\noindent{\bf Proof.} The proof of inequality~\eqref{Fsumruleif} is similar to that of~\eqref{FC2}. To prove~\eqref{FsumruleifE} under each of the conditions (a)--(b), define the sets
	\begin{align*}
		\Omega_1=\big\{(x, y_1, y_2)\in X\times Y\times Y \mid y_1\in F_1(x)\big\},\\
		\Omega_2=\big\{(x, y_1, y_2)\in X\times Y\times Y \mid y_2\in F_2(x)\big\}.
	\end{align*}
	(a) Under the assumptions made, we have $\mbox{\rm int}(\Omega_1)\cap \Omega_2\neq\emptyset$. Indeed, if $\mbox{\rm int}(\Omega_1)\cap \Omega_2=\emptyset$, then by convex separation there exists a nonzero element $(x^*, y_1^*, y^*_2)$ such that
	\begin{equation*} \label{separation_ineq}
		\la x^*, x\ra+\la y_1^*, y_1\ra+\la y^*_2, y_2\ra \leq  \la x^*, u\ra+\la y_1^*, z_1\ra+\la y^*_2, z_2\ra
	\end{equation*}
	whenever $(x, y_1, y_2)\in \Omega_1$ and $(u, z_1, z_2)\in \Omega_2$. The constructions of $\Omega_1$ and $\Omega_2$ yield $y^*_1=y^*_2=0$. Then $x^*\neq 0$ and
	\begin{equation}\label{eqs1}
		\la x^*, x\ra\leq \la x^*, u\ra\ \; \mbox{\rm whenever }x\in \dom F_1,\; u\in \dom F_2.
	\end{equation}
	Choose $\ox\in \mbox{\rm int}(\dom F_1)\cap \dom F_2$. Then we find $\delta>0$  and a neighborhood $U$ of the origin in $X$ such that $\ox+\delta U\subset \dom F_1$. It follows from~\eqref{eqs1} that
	\begin{equation*}
		\la x^*, \ox+\delta b\ra\leq \la x^*, \ox\ra\ \; \mbox{\rm whenever }\; b\in  U,
	\end{equation*}
	which implies that $\la x^*, b\ra\leq 0$ whenever $b\in U$. This clearly gives $x^*=0$, so we have $(x^*, y^*_1, y^*_2)=(0,0,0)$, a contradiction.  Following the proof of Theorem~\ref{Fsumrule}, it suffices to establish the quality~\eqref{FsumruleifE} for the case where $(u^*, v^*)\in \mbox{\rm dom}(F_1+F_2)^*$. Arguing similarly as in the proof of Lemma~\ref{lem:conj_sum}, we can show that
	\begin{equation}\label{equality_Fenchel}
		(F_1+F_2)^*(u^*, v^*)=\sigma_{\small {\Omega_1\cap \Omega_2}}(u^*, v^*).
	\end{equation}
	Now, we only need to repeat the proof of Theorem~\ref{Fsumrule} using Theorem~\ref{Flemma1}(a), where the condition ~$\Omega_1\cap \mbox{\rm int}(\Omega_2)\neq\emptyset$ is replaced with~$\mbox{\rm int}(\Omega_1)\cap \Omega_2\neq\emptyset$. \\[1ex]
	(b) Since $\gph F_2$ is a polyhedral convex set, we see that $\Omega_2$ is a polyhedral convex set. By the assumptions, $\hat{x}\in\mbox{\rm qri}(\dom F_1)\cap \dom F_2$ and $(\hat{x}, \hat{y}_1)\in \mbox{\rm qri}(\gph F_1)$. Choose $\hat{y}_2\in F_2(\hat{x})$. By the constructions of $\Omega_1$ and $\Omega_2$ we have $(\hat{x}, \hat{y}_1, \hat{y}_2)\in \mbox{\rm qri}(\Omega_1)\cap \Omega_2$, so $\mbox{\rm qri}(\Omega_1)\cap \Omega_2\neq\emptyset$.  Similar to part (a), we can employ Theorem~\ref{Flemma1}(b), where the roles of $\Omega_1$ and $\Omega_2$ are swapped, to finish the proof. \\[1ex]
	(c) Consider the case where both $F_1$ and $F_2$ are polyhedral convex set-valued mappings with $\dom F_1\cap \dom F_2\neq\emptyset$. Then $\Omega_1$ and $\Omega_2$ are polyhedral convex sets. Fix $\hat{x}\in \dom F_1\cap \dom F_2$. Then choose $\hat{y}_1\in F_1(\hat{x})$ and $\hat{y}_2\in F_2(\hat{x})$. It follows easily that $(\hat{x}, \hat{y}_1, \hat{y}_2)\in \Omega_1\cap \Omega_2$, so $\Omega_1\cap \Omega_2\neq\emptyset$. Now we can finish the proof using Theorem~\ref{Flemma1}(c).\\[1ex]
	(d) We have the equality
	\begin{equation}\label{difference}
		\Omega_1-\Omega_2=(\dom F_1-\dom F_2)\times Y\times Y.
	\end{equation}
	Indeed, given any $(x,y,z)\in (\dom F_1-\dom F_2)\times Y\times Y$, we can find $x_1\in \dom F_1$ and $x_2\in\dom F_2$ such that $x=x_1-x_2$. Take any $y_1\in F_1(x_1)$, $z_2\in F_2(x_2)$, and put $z_1=y_1-y$, $y_2= z_2-z$. Since $(x_1,y_1,y_2)\in \Omega_1$, $(x_2,z_1,z_2)\in\Omega_2$, and $(x,y,z)=(x_1,y_1,y_2)- (x_2,z_1,z_2)$, the inclusion $\supset$ is valid. To prove the reverse inclusion, take any $(x_1,y_1,y_2)\in \Omega_1$, $(x_2,z_1,z_2)\in\Omega_2$, and note that $x_1\in \dom F_1$ and $x_2\in\dom F_2$. Therefore,
	$$(x_1,y_1,y_2)- (x_2,z_1,z_2) =(x_1-x_2,y_1-z_1,y_2-z_2)\in (\dom F_1-\dom F_2)\times Y\times Y.$$
	It follows that $\Omega_1-\Omega_2\subset(\dom F_1-\dom F_2)\times Y\times Y$. As $F_1$ and $F_2$ have closed convex graphs, $\Omega_1$ and $\Omega_2$ are closed convex sets. By~\eqref{regu_d}, $\cone(\dom F_1-\dom F_2)$  is a closed subspace. Observe that
	\begin{equation*}
		\cone(\Omega_1-\Omega_2)=\cone(\dom F_1-\dom F_2)\times Y\times Y.
	\end{equation*}
	Indeed, the inclusion $\subset$ is obvious. Now, take any $$(x, y_1, y_2)\in \cone(\dom F_1-\dom F_2)\times Y\times Y.$$ Since $\cone(\dom F_1-\dom F_2)$ is a closed linear subspace, by \cite[Lemma~1.4]{AB} we have
	\begin{equation*}
		\cone(\dom F_1-\dom F_2)=\bigcup_{t>0}t\big(\dom F_1-\dom F_2\big).
	\end{equation*}
	So, there exist $t>0$ and $a\in \dom F_1-\dom F_2$ such that $x=ta$. Thus, thanks to~\eqref{difference}, we have
	\begin{equation*}
		(x, y_1, y_2)=t(a, y_1/t, y_2/t)\in \cone(\Omega_1-\Omega_2).
	\end{equation*}
	Combining~\eqref{difference} with~\eqref{regu_d} yields $0\in \mbox{\rm sqri}(\Omega_1-\Omega_2)$. Following the proof of Theorem~\ref{Fsumrule}, it suffices to prove the quality~\eqref{FsumruleifE} for the case where $(u^*, v^*)\in \mbox{\rm dom}(F_1+F_2)^*$. Since~\eqref{equality_Fenchel} is valid, it remains to repeat the proof of Theorem~\ref{Fsumrule} using Theorem~\ref{Flemma1}(d). $\h$

We now establish the Fenchel conjugate chain rule for set-valued mappings.
\begin{theorem}\label{FchainruleIN}
	Let $F\colon X\tto Y$ and $G\colon Y\tto Z$ be set-valued mappings, and let $(x^*,z^*)\in X^*\times Z^*$. Then we have
	\begin{equation}\label{FchainruleIF}
		(G\circ F)^*(x^*,z^*) \leq \inf\{F^*(x^*,v^*)+G^*(-v^*,z^*)\mid v^*\in Y^*\}.
	\end{equation}
	Suppose in addition that one of the following conditions is satisfied:
	\begin{enumerate}
		\item $F$ and $G$ are convex set-valued mappings such that
		\begin{equation}\label{chain_rule_2}
			\mbox{\rm int}(\gph F)\neq\emptyset \ \; \mbox{\rm and } \ \;  \mbox{\rm int}(\rge F)\cap (\dom G)\neq\emptyset.
		\end{equation}
		\item $F$ is a polyhedral convex set-valued mapping and $G$ is a convex set-valued mapping such that there exist $\hat{y}\in \rge F$ and $\hat{z}\in Z$ with $(\hat{y}, \hat{z})\in \mbox{\rm qri}(\gph G)$ (or $F$ is a convex set-valued mapping and $G$ is a polyhedral convex set-valued mapping such that  there exist $\hat{x}\in X$ and  $\hat{y}\in  \dom G$  such that $(\hat{x}, \hat{y})\in \mbox{\rm qri}(\gph F)$).
		\item Both $F$ and $G$ are polyhedral convex set-valued mappings and $\rge F\cap \dom G\neq\emptyset$.
		\item $X$ and $Y$ are Banach spaces, $F$ and $G$ are convex set-valued mappings with  closed graphs, and $0\in \sqri(\rge F-\dom G).$
	\end{enumerate}
	Then we have the equality
	\begin{equation}\label{FchainruleIFE}
		(G\circ F)^*(x^*,z^*) = \inf\{F^*(x^*,v^*)+G^*(-v^*,z^*)\mid v^*\in Y^*\}.
	\end{equation}
	Furthermore, the infimum in~\eqref{FchainruleIFE} is attained if $(x^*, z^*)\in \mbox{\rm dom}(G\circ F)^*$.
\end{theorem}
\noindent{\bf Proof.} To prove~\eqref{FchainruleIF}, it suffices to repeat the proof of~\eqref{Fchainruleff} from Theorem~\ref{Fchainrulef}. \\[1ex]
	Now consider the case where $F$ and $G$ are convex. Define two convex sets
	\begin{align*}
		&\Omega_1=\big\{(x, y, z)\; \big |\; y\in F(x)\big\}=(\gph F)\times Z,\\
		&\Omega_2=\big\{(x, y, z)\; \big |\; z\in F(y)\big\}=X\times \gph G.
	\end{align*}
	(a) Assume that	the two conditions in~\eqref{chain_rule_2} are satisfied. Let us show that $\mbox{\rm int}(\Omega_1)\cap \Omega_2\neq\emptyset$. Suppose on the contrary that
	\begin{equation}\label{intersection_int}
		\mbox{\rm int}(\Omega_1)\cap \Omega_2=\emptyset.\end{equation}
	Since $\mbox{\rm int}(\Omega_1)\neq\emptyset$, by convex separation, there is a nonzero triple $(x^*, y^*, z^*)\in X^*\times Y^*\times Z^*$ such that
	\begin{equation*}
		\la x^*, x_1\ra+\la y^*, y_1\ra+\la z^*, z_1\ra\leq \la x^*, x_2\ra+\la y^*, y_2\ra+\la z^*, z_2\ra
	\end{equation*}
	whenever $(x_1, y_1, z_1)\in \Omega_1$ and 	$(x_2, y_2, z_2)\in \Omega_2$. By the constructions of $\Omega_1$ and $\Omega_2$ we have $x^*=0$ and $z^*=0$. Then
	\begin{equation*}
		\la  y^*, y_1\ra\leq \la y^*, y_2\ra\ \; \mbox{\rm whenever }\; y_1\in \rge F,\ y_2\in \dom G.
	\end{equation*}
	Since $\mbox{\rm int}(\rge F)\cap \dom G\neq\emptyset$, we can deduce that $y^*=0$, which is a contradiction because $(x^*, y^*, z^*)\neq 0$. Now, we can repeat the second part of the proof of Theorem~\ref{Fchainrulef}, where~\eqref{ri_intersection} is replaced by~\eqref{intersection_int}. With the use of Theorem~\ref{Flemma1} under~(a), where the roles of $\Omega_1$ and $\Omega_2$ are swapped, we can derive the conclusions.\\[1ex]
	(b) Under the assumptions made, choose $\hat{y}\in \rge F$ and $\hat{z}\in Z$ with $(\hat{y}, \hat{z})\in \mbox{\rm qri}(\gph G)$. Then choose $\hat{x}\in X$ such that $(\hat{x}, \hat{y})\in \gph F$. We have $(\hat{x}, \hat{y}, \hat{z})\in \Omega_1\cap \mbox{\rm qri}(\Omega_2).$
	Since $F$ is polyhedral convex, so is the set $\Omega_1$. Now, we can complete the proof in this case using the arguments of the proof of Theorem~\ref{Fchainrulef} and Theorem~\ref{Flemma1} under (b). \\[1ex]
	(c) By the assumptions made, we can easily show that $\Omega_1$ and $\Omega_2$ are polyhedral convex sets with nonempty intersection. It remains to use the method for proving Theorem~\ref{Fchainrulef} and Theorem~\ref{Flemma1} under (c) to finish the proof in this case. \\[1ex]
	(d) Under our assumptions, $\Omega_1$ and $\Omega_2$ are nonempty closed convex sets in the Banach space $X\times Y\times Z$. Note that
	\begin{equation*}\label{rgeF-domG}
		\Omega_1-\Omega_2=X\times (\rge F-\dom G)\times Z
	\end{equation*}
	and
	\begin{equation}\label{cone f}
		\cone(\Omega_1-\Omega_2)=X\times \cone(\rge F-\dom G)\times Z.
	\end{equation}
	Since $0\in \sqri(\rge F-\dom G)$, $\cone(\rge F-\dom G)$ is closed linear subspace of $Y$. Then, from~\eqref{cone f} it follows that $\cone(\Omega_1-\Omega_2)$ is a closed linear subspace of $X\times Y\times Z$. This means that $0\in\mbox{\rm sqri}\,(\Omega_1-\Omega_2)$. So,
	we can complete the proof in this case using an analogue of Lemma~\ref{lem:chainrule} and Theorem~\ref{Flemma1} under~(d).$\h$

\medskip
The last result of this section concerns the Fenchel conjugate intersection rule for set-valued mappings in infinite dimensions.

\begin{theorem}
	Let $F_1, F_2\colon X\tto Y$ be two set-valued mappings. Then
	\begin{equation*}
		(F_1\cap F_2)^*(x^*, y^*)\leq (F_1^*\square F_2^*)(x^*, y^*)\ \; \mbox{\rm for all }\; (x^*, y^*)\in X^*\times Y^*.
	\end{equation*}
	Assume that one of the following conditions is satisfied:
	\begin{enumerate}
		\item $\mbox{\rm int}(\gph F_1)\cap \gph F_2\neq\emptyset$.
		\item $F_2$ is polyhedral convex and $\mbox{\rm qri}(\gph F_1)\cap \gph F_2\neq\emptyset$.
		\item $F_1$ and $F_2$ are polyhedral convex and $\gph F_1\cap \gph F_2\neq\emptyset$.
		\item $X$ and $Y$ are Banach spaces,  $0\in \sqri(\gph F_1-\gph F_2)$, $F_1$ and $F_2$ have closed graphs.
	\end{enumerate}\label{FInruleE}
	Then we have the equality
	\begin{equation*}
		(F_1\cap F_2)^*(x^*, y^*)= (F_1^*\square F_2^*)(x^*, y^*)\ \; \mbox{\rm for all }\; (x^*, y^*)\in X^*\times Y^*
	\end{equation*}
	Furthermore, if $(F_1\cap F_2)^*(x^*, y^*)\in \R$, then the infimum in the infimal convolution from~\eqref{FInruleE} is attained.
\end{theorem}
\noindent{\bf Proof.} We have by definition that
	\begin{equation*}
		(F_1\cap F_2)^*(x^*,y^*) = 	\sigma_{\mbox{\rm gph}(F_1\cap F_2)}(x^*,y^*) =\sigma_{\gph F_1\cap\gph F_2}(x^*,y^*)
	\end{equation*} for all $(x^*, y^*)\in X^*\times Y^*$. To complete the proof, we only need to apply Theorem~\ref{Flemma1}.
$\h$

\section{Derivation of Convex Generalized Differentiation}
\setcounter{equation}{0}

In this section, we will show that major known calculus rules of convex generalized differentiation can be derived from calculus rules of Fenchel conjugate for set-valued mappings.

In the theorem below, we derive the coderivative sum rule for set-valued mappings using the corresponding Fenchel conjugate sum rule. We refer the reader to~\cite{bmn,bmn2022} for several coderivative sum rules which were established by other proofs.

\begin{theorem}\label{coderivative_sum_rule} Let $F_1\colon X\tto Y$ and $F_2\colon X\tto Y$ be convex set-valued mappings. Suppose that one of the following conditions is satisfied:
	\begin{enumerate}
		\item $X=\R^n$, $Y=\R^p$, and  one of the assumptions (a), (b), or (c) from Theorem~\ref{Fsumrule} is fulfilled.
		\item One of the assumptions (a), (b), (c),  or (d) from Theorem~\ref{Fsumrulei} is fulfilled.
	\end{enumerate}
	Then, for any $(\ox, \oy)\in \mbox{\rm gph}(F_1+F_2)$ and $(\oy_1, \oy_2)\in S(\ox, \oy)$ with \begin{equation*}
		S(\ox, \oy):=\big\{(\oy_1, \oy_2)\in Y\times Y\; \big |\; \oy_i\in F_i(\ox), \; \oy=\oy_1+\oy_2\big\},
	\end{equation*} we have the equality
	\begin{equation}\label{Csumrule}
		D^*(F_1+F_2)(\ox, \oy)(y^*)=D^*F_1(\ox, \oy_1)(y^*)+D^*F_2(\ox, \oy_2)(y^*)
	\end{equation}
	for every $y^*\in Y^*$.
\end{theorem}
\noindent{\bf Proof.} It suffices to prove the inclusion $\subset$ in~\eqref{Csumrule} because the reverse inclusion can be easily verified by using the definitions of coderivative  in~\eqref{def_coderivative} and normal cone in~\eqref{def_normal_cone}.
	
	Fix any $x^*\in  D^*(F_1+F_2)(\ox, \oy)(y^*)$. Then by Theorem~\ref{T1} we have
	\begin{equation}\label{YIG}
		\la x^*, \ox\ra= \la y^*, \oy\ra +(F_1+F_2)^*(x^*, -y^*).
	\end{equation}
	Since the last equality implies that $(F_1+F_2)^*(x^*, -y^*)\in \R$, we can apply Theorem~\ref{Fsumrule} or Theorem~\ref{Fsumrulei} and find $x^*_1, x^*_2\in X^*$ with $x^*=x^*_1+x^*_2$ and
	\begin{equation*}
		(F_1+F_2)^*(x^*, -y^*)=F_1^*(x^*_1, -y^*)+F^*_2(x^*_2, -y^*).
	\end{equation*}
	It follows from~\eqref{YIG} that
	\begin{equation*}
		\la x^*_1, \ox\ra+\la x^*_2, \ox\ra= \la y^*, \oy_1\ra+ \la y^*, \oy_2\ra +F_1^*(x^*_1, -y^*)+F^*_2(x^*_2, -y^*).
	\end{equation*}
	We can deduce from this using~\eqref{Young1} that
	\begin{equation*}
		\la x^*_1, \ox\ra=\la y^*, \oy_1\ra+F_1^*(x^*_1, -y^*),\ \;
		\la x^*_2, \ox\ra=\la y^*, \oy_2\ra+F_2^*(x^*_2, -y^*).
	\end{equation*}
	Applying Theorem~\ref{T1} again, we have
	\begin{equation*}
		x^*_1\in D^*F_1(\ox, \oy_1)(y^*),\ \; x^*_2\in D^*F_2(\ox, \oy_2)(y^*).
	\end{equation*}
	Therefore,
	\begin{equation*}
		x^*=x^*_1+x^*_2\in D^*F_1(\ox, \oy_1)(y^*)+D^*F_2(\ox, \oy_2)(y^*),
	\end{equation*}
	which completes the proof of the theorem. $\h$

Theorem~\ref{coderivative_sum_rule} can be used to obtain the Fenchel conjugate sum rule and the subdifferential sum rule for extended-real-valued functions.

\begin{theorem} Let $f_1, f_2\colon X\to \oR$ be two proper convex functions. Suppose that one of the following conditions is satisfied:
	\begin{enumerate}
		\item $X=\R^n$ and $\mbox{\rm ri}(\dom f_1)\cap \mbox{\rm ri}(\dom f_2)\neq\emptyset$.
		\item $X=\R^n$, $f_2$ is polyhedral convex and $\mbox{\rm ri}(\dom f_1)\cap \dom f_2\neq\emptyset$.
		\item $f_1$ is continuous at some point $\ou\in \dom f_1\cap \dom f_2$.
		\item $f_2$ is polyhedral convex and $\mbox{\rm qri}(\dom f_1)\cap \dom f_2\neq\emptyset$.
		\item $f_1$ and $f_2$ are polyhedral convex and $\dom f_1\cap \dom f_2\neq\emptyset$.
		\item $X$ is a Banach space,  $f_1$ and $f_2$ are lower semicontinuous, and $$0\in \sqri(\dom f_1-\dom f_2).$$
	\end{enumerate}
	Then, for every $x^*\in X^*$, we have the equality
	\begin{equation*}\label{FsumruleFunction}
		(f_1+f_2)^*(x^*)=(f^*_1\s f^*_2)(x^*),
	\end{equation*}
	where the infimum in \begin{equation}\label{equality_n1} (f^*_1\s f^*_2)(x^*)=\inf\big\{f^*_1(x^*_1)+f^*_2(x^*_2)\; \big|\; x^*_1+x^*_2=x^*\big\}
	\end{equation}
	is attained if $(f_1+f_2)^*(x^*)\in \R$. In addition,
	\begin{equation}\label{SsumruleFunction}
		\partial (f_1+f_2)(\ox)=\partial f_1(\ox)+\partial f_2(\ox)
	\end{equation} for  any $\ox\in \dom f_1\cap \dom f_2$.
\end{theorem}
\noindent{\bf Proof.} Consider the epigraphical mappings
	\begin{equation*}
		E_{f_i}(x)=\big\{\lambda \in \R\; \big |\; f_i(x)\leq \lambda\big\},\ \; x\in X,
	\end{equation*}
	for $i=1,2$. Then $\gph E_{f_i}=\epi f_i$ and $\dom E_{f_i}=\dom f_i$. Let $F_i=E_{f_i}$ for $i=1, 2$.
	Suppose that~(a) or~(b) is satisfied. Then it is obvious that either assumption~(a) or assumption~(b) in Theorem~\ref{Fsumrule} is satisfied. Now, suppose that~(c) is satisfied. Then $\mbox{\rm int}(\epi f_1)\neq\emptyset$, and $\ou \in \mbox{\rm int}(\dom f_1)\cap \dom f_2$. Thus, $\mbox{\rm int}(\gph E_{f_1})\neq\emptyset$ and $\mbox{\rm int}(\dom E_{f_1})\cap \dom E_{f_2}\neq\emptyset$, so assumption~(a) of Theorem~\ref{Fsumrulei} is satisfied. Next, assume that~(d) is satisfied. Fix $\hat{x}\in \mbox{\rm qri}(\dom f_1)\cap \dom f_2$ and choose a real number $\hat{\lambda}$ such that  $\hat{\lambda}>\max\{f_1(\hat{x}), f_2(\hat{x})\}$. Then by~\cite[Theorem~2.190]{bmn2022} one has $(\hat{x}, \hat{\lambda})\in \mbox{\rm qri}(\epi f_1)$. Thus, assumption~(b) in Theorem~\ref{Fsumrulei} is satisfied. Under~(e) we can easily see that assumption~(c) in Theorem~\ref{Fsumrulei} is satisfied. Finally, if~(f) is satisfied, then assumption~(d) in Theorem~\ref{Fsumrulei} is satisfied. Therefore, by Proposition~\ref{fF}, Theorem~\ref{Fsumrule}, or Theorem~ \ref{Fsumrulei} we have
	\begin{align*}
		(f_1+f_2)^*(x^*)&=(E_{f_1}+E_{f_2})^*(x^*, -1)\\
		& =\inf\big\{(E_{f_1})^*(x^*_1, -1)+(E_{f_2})^*(x^*_2, -1)\; \big |\;  x^*_1+x^*_2=x^*\big\}\\
		&=\inf\big\{f^*_1(x^*_1)+f^*_2(x^*_2)\; \big |\;  x^*_1+x^*_2=x^*\big\}\\
		&=(f^*_1\s f^*_2)(x^*).
	\end{align*}
	Furthermore, the infimum in~\eqref{equality_n1} is attained if $(f_1+f_2)^*(x^*)\in \R$ because, in that case, $(x^*, -1)\in \mbox{\rm dom}(E_{f_1}+E_{f_2})^*$ and the first infimum in the last transformation is achieved either by Theorem~\ref{Fsumrule} or Theorem~ \ref{Fsumrulei}.
	
	Now take any $\ox\in \dom f_1\cap \dom f_2$. Since the functions $f_1$ and $f_2$ are proper by our assumptions, a direct verification shows that
	$$E_{f_1+f_2}(x)=E_{f_1}(x)+E_{f_2}(x)$$
	for all $x\in X$. So, applying formula~\eqref{coder_sub} gives
	\begin{equation*}
		\begin{aligned}
			\partial (f_1+f_2) (\ox) &=D^*(E_{f_1+f_2})(\ox, (f_1+f_2)(\ox))(1)\\
			&=D^*(E_{f_1}+E_{f_2})(\ox,f_1(\ox)+f_2(\ox))(1).
		\end{aligned}
	\end{equation*}
	Then, under each of the assumptions~(a)--(f), we can apply Theorem~\ref{coderivative_sum_rule} and get from formula~\eqref{Csumrule} that
	\begin{equation*}
		\begin{aligned}
			\partial (f_1+f_2) (\ox) &= D^*(E_{f_1})(\ox, f_1(\ox))(1)+  D^*(E_{f_2})(\ox, f_2(\ox))(1)\\
			&=\partial f_1(\ox)+\partial f_2(\ox).
		\end{aligned}
	\end{equation*}
	We have thus proved the equality~\eqref{SsumruleFunction}. $\h$

Given two set-valued mappings $F\colon X\tto Y$ and $G\colon Y\tto Z$, define
\begin{equation*}
	M(x, z)=F(x)\cap G^{-1}(z),\ \; (x, z)\in X\times Z.
\end{equation*}

In the theorem below, we derive the coderivative chain rule for convex set-valued mappings from the corresponding Fenchel conjugate chain rule.
\begin{theorem}\label{F_chain_rule} Let $F\colon X\tto Y$ and $G\colon Y\tto Z$ be convex set-valued mappings. Suppose that one of the following conditions is satisfied:
	\begin{enumerate}
		\item $X=\R^n$, $Y=\R^p$, $Z=\R^q$, and  one of the assumptions (a), (b), or (c) from Theorem~\ref{Fchainrulef} is fulfilled.
		\item $X$,  $Y$, $Z$ are LCTV spaces and one of the assumptions (a), (b), (c),  or (d) from Theorem~\ref{FchainruleIN} is fulfilled.
	\end{enumerate}
	Then for any $(\ox, \oz)\in \mbox{\rm gph}(G\circ F)$ and $z^*\in Z^*$, we have the equality
	\begin{equation}\label{Cchainrule}
		D^*(G\circ F)(\ox, \oz)(z^*)=\big(D^*F(\ox, \oy)\circ D^*G(\oy, \oz)\big)(z^*)
	\end{equation}
	whenever $\oy \in M(\ox, \oz)$.
\end{theorem}
\noindent{\bf Proof.} We first verify the inclusion $\subset$ in~\eqref{Cchainrule}.  Since $F$ and $G$ are convex set-valued mappings, a direct verification shows that the composition $G\circ F\colon X\tto Z$ is a convex set-valued mapping. Fix any $x^*\in D^*(G\circ F)(\ox, \oz)(z^*)$ and any $\oy\in M(\ox, \oz)$.  Then by Theorem~\ref{T1} we have
	\begin{equation*}\label{YIG_a}
		\la x^*, \ox\ra= \la z^*, \oz\ra +(G\circ F)^*(x^*, -z^*).
	\end{equation*}
	This yields $(G\circ F)^*(x^*, -z^*)\in \R$. So, by the assumptions made and by Theorem~\ref{Fchainrulef} or Theorem~\ref{FchainruleIN} we have
	\begin{equation*}
		\la x^*, \ox\ra= \la z^*, \oz\ra +F^*(x^*, -v^*)+G^*(v^*, -z^*)
	\end{equation*}
	for some $v^*\in Y^*$.  It follows that
	\begin{equation*}
		\la x^*, \ox\ra + \la v^*, \oy\ra= \la v^*, \oy\ra+ \la z^*, \oz\ra +F^*(x^*, -v^*)+G^*(v^*, -z^*).
	\end{equation*}
	We can deduce from this using~\eqref{Young1} that
	\begin{align*}
		&\la x^*, \ox\ra=\la v^*, \oy\ra+F^*(x^*, -v^*),\\
		&\la v^*, \oy\ra=\la z^*, \oz\ra+G^*(v^*, -z^*).
	\end{align*}
	As $F$ and $G$ are convex set-valued mappings, applying Theorem~\ref{T1} again gives
	\begin{equation*}
		x^*\in D^*F(\ox, \oy)(v^*),\ \; v^*\in D^*G(\oy, \oz)(z^*).
	\end{equation*}
	Therefore, we have by definition that
	\begin{equation*}
		x^*\in \big(D^*F(\ox, \oy)\circ D^*G(\oy, \oz)\big)(z^*),
	\end{equation*}
	which implies that $D^*(G\circ F)(\ox, \oz)(z^*)\subset\big(D^*F(\ox, \oy)\circ D^*G(\oy, \oz)\big)(z^*)$.
	
	To obtain the reverse inclusion and, therefore, complete the proof of~\eqref{Cchainrule}, take any $x^*\in \big(D^*F(\ox, \oy)\circ D^*G(\oy, \oz)\big)(z^*)$. Let  $y^*\in D^*G(\oy, \oz)(z^*)$ be a vector such that $x^*\in D^*F(\ox, \oy)(y^*)$. Then, as $F$ and $G$ are convex set-valued mappings, using~\eqref{def_coderivative} yields
	\begin{equation*}
		\begin{cases}
			(x^*,-y^*)\in N\big((\ox,\oy);\gph F\big),\\
			(y^*,-z^*)\in N\big((\oy,\oz);\gph G\big)
		\end{cases}
	\end{equation*}
	or, equivalently,
	\begin{equation}
		\label{eqn:chainrule}
		\begin{cases}
			\la x^*,x-\ox \ra -\la y^*,y-\oy\ra \leq 0,\;\; \mbox{for all}\; (x,y)\in \gph F,\\
			\la y^*, y-\oy\ra - \la z^*,z-\oz\ra \leq 0,\;\; \mbox{for all}\; (y,z)\in \gph G.
		\end{cases}
	\end{equation}
	For every $(x,z)\in \mbox{\rm gph}(G\circ F)$, there is some $y\in F(x)$ such that $z\in G(y)$. Hence, adding the inequalities in~\eqref{eqn:chainrule} side-by-side gives $$\la x^*,x-\ox \ra + \la -z^*,z-\oz\ra \leq 0.$$
	Since $(x,z)\in \mbox{\rm gph}(G\circ F)$ can be chosen arbitrarily, this inequality indicates that $$(x^*,-z^*)\in N\big((\ox,\oz); \mbox{\rm gph}(G\circ F)\big).$$ Therefore, $x^*\in D^*(G\circ F)(\ox, \oz)(z^*)$.$\h$

From Theorem~\ref{Fchainrulef}, Theorem~\ref{FchainruleIN}, and Theorem~\ref{F_chain_rule} we can obtain the next Fenchel conjugate chain rule and subdifferential chain rule for extended-real-valued functions.

\begin{theorem} Let $A\colon X\to Y$ be a continuous linear operator, and let  $g\colon Y\to \oR$ be a convex function. Suppose that one of the following conditions is satisfied:
	\begin{enumerate}
		\item $X=\R^n$, $Y=\R^p$ and $A(X)\cap \mbox{\rm ri}(\dom g)\neq\emptyset$.
		\item $g$ is continuous at some point $\ou\in A(X)$.\
		\item $g$ is polyhedral convex and $A(X)\cap \dom g\neq\emptyset$.
		\item $X,Y$ are Banach spaces, $g$ is lower semicontinuous, and $0\in \sqri(A(X)-\dom g)$.
	\end{enumerate}
	Then we have the equality
	\begin{equation}\label{FchainruleFunction}
		(g\circ A)^*(x^*)=\inf\big\{g^*(y^*)\; \big| \; A^*y^*=x^*\big\},
	\end{equation}
	where the infimum in~\eqref{FchainruleFunction} is attained if $(g\circ A)^*(x^*)\in \R$. In addition,
	\begin{equation}\label{SchainruleFunction}
		\partial (g\circ A)(\ox)=A^*\big(\partial g(A\ox)\big),\ \; \ox\in \mbox{\rm dom}(g\circ A).
	\end{equation}
\end{theorem}
\noindent{\bf Proof.} Put $F(x)=\{Ax\}$ for all $x\in X$, $G(y)=E_g(y)$ for all $y\in Y$. Clearly, $\rge F=A(X)$ and by~\eqref{Emapping} one has $\dom G=\dom g$. Furthermore, $\mbox{\rm epi}(g\circ A)=\mbox{\rm gph}(G\circ F)$. Since the assumptions in either Theorem~\ref{Fchainrulef}, or Theorem~\ref{FchainruleIN}, and Theorem~\ref{F_chain_rule} are satisfied under our setting, by Proposition~\ref{fF} and by one of these theorems we have
	\begin{equation*}\begin{array}{rl} (g\circ A)^*(x^*) =E^*_{g\circ A}(x^*, -1) & =(G\circ F)^*(x^*, -1)\\
			& =\inf\big\{ F^*(x^*, -v^*)+G^*(v^*, -1)\; \big |\; v^*\in Y^*\big\},		
		\end{array}
	\end{equation*} where the infimum is attained if $(g\circ A)^*(x^*)\in \R$. As $F^*(x^*, -v^*)=\infty$ if $A^*v^*\neq x^*$  and $F^*(x^*, -v^*)=0$ if $A^*v^*=x^*$ (see Example~\ref{Flinear}), the latter implies that
	\begin{align*}
		(g\circ A)^*(x^*)&=\inf\big\{ F^*(x^*, -v^*)+G^*(v^*, -1)\; \big |\; v^*\in Y^*\big\}\\
		&=\inf\big\{G^*(y^*, -1)\; \big |\; A^*y^*=x^*\big\}\\
		&= \inf\big\{g^*(y^*)\; \big| \; A^*y^*=x^*\big\},
	\end{align*}  where the last infimum is attained if $(g\circ A)^*(x^*)\in \R$. Thus, the first assertion of the theorem is  true.
	
	For any $\ox\in \mbox{\rm dom}(g\circ A)$, by~\eqref{coder_sub} and the equality $\mbox{\rm epi}(g\circ A)=\mbox{\rm gph}(G\circ F)$ we have
	$$\partial (g\circ A)(\ox)=D^*E_{g\circ A}(\ox,(g\circ A)(\ox))(1)=D^*(G\circ F)(\ox,(g\circ A)(\ox))(1).$$ So, applying Theorem~\ref{F_chain_rule} and formula~\eqref{coder_sub} gives
	$$\begin{array}{rl}
		\partial (g\circ A)(\ox)&=\big(D^*F(\ox,A\ox)\circ D^*G(A\ox,g(A\ox))\big)(1)\\
		&=D^*F(\ox,A\ox)\big(D^*E_g(A\ox,g(A\ox))(1)\big)\\
		&=D^*F(\ox,A\ox)\big(\partial g(A\ox)\big)\\
		&=A^*\big(\partial g(A\ox)\big).
	\end{array}$$
	This establishes~\eqref{SchainruleFunction} and completes the proof of the theorem. $\h$

\section{Concluding Remarks}
\label{s:Concluding}

We have introduced a new notion of Fenchel conjugate for set-valued mappings in infinite dimensions and provided a comprehensive study of its properties and calculus rules. The study highlights the importance of qualification conditions involving relative interiors, interiors, and their generalizations, while also demonstrating that these conditions can be relaxed with the presence of polyhedral convexity.

The obtained results have significant implications in convex analysis, as they allow the unified  derivation of old and new calculus rules for coderivatives of set-valued mappings as well as Fenchel  conjugates and subdifferentials of extended-real-valued functions.  It is likely that this approach also enables us to get a refined coderivative rule for the intersection of two set-valued mappings, along with enhanced Fenchel conjugate and subdifferential formulas for the maximum of convex functions.

This research contributes to advancing our understanding of set-valued mapping theory in infinite dimensions and opens up avenues for further exploration in the area.

\end{document}